\documentclass[a4, 12pt]{article}
\usepackage{amssymb}
\usepackage{yfonts}
\def\pointir{\unskip . --- \ignorespaces}

\newtheorem{proposition}{Proposition}[section]
\newtheorem{corollary}{Corollary}[section]
\newtheorem{lemma}{Lemma}[section]
\newtheorem{theorem}{Theorem}[section]

\makeatletter

\title{\bf Multiplicity one Theorems}
\author{Avraham Aizenbud,
\and Dmitry Gourevitch,
\and Steve Rallis \footnote{The third author is partially supported by NSF Grant DMS-0500392},
\and G\'erard Schiffmann}

\begin{document}

\def\section
{\@startsection {section}{1}
{-1pt}{-5ex \@plus -1ex \@minus -.2ex}
{2ex \@plus .2ex}
{\normalfont \large \bfseries } }

\maketitle

\vskip 0.5cm
 \begin{abstract}
In the local, characteristic 0, non archimedean case,  we consider
distributions on $GL(n+1)$ which are invariant under the adjoint
action of  $GL(n)$. We prove that such distributions are invariant
by transposition. This  implies that an admissible irreducible representation of
$GL(n+1)$, when restricted to $GL(n)$ decomposes with multiplicity one.\hfill\break
 Similar Theorems are obtained for orthogonal or unitary
groups.\end{abstract}

 \centerline{\bf Introduction}
 \vskip 0.5cm
Let ${\mathbb F }$ be a local field non archimedean and of
characteristic 0. Let $W$ be a vector space over ${\mathbb F}$ of
finite dimension $n+1\geqslant  1$ and let $W=V\oplus U$ be a
direct sum decomposition with $\dim V=n $. Then we have an
imbedding of $GL(V)$ into $GL(W)$. Our goal is to prove the
following Theorem:

\vskip 0.5cm
{\parindent 0pt
{\bf Theorem $\mathbf 1$:}
 {\sl If $\pi$ (resp. $\rho$) is an irreducible admissible representation of $GL(W)$ (resp. of $GL(V)$) then
\[ \mathrm{dim} \left ({\mathrm {Hom}}_ {GL(V)}(\pi _{|GL(V)}\, ,\, \rho )\right )\leqslant 1.\] }}

We choose a basis of $V$ and a non zero vector in $U$ thus getting a basis of $W$. We can identify $GL(W)$ with $GL(n+1,{\mathbb F})$ and $GL(V)$ with $GL(n,{\mathbb F})$. The transposition map is an involutive anti-automorphism of $GL(n+1,{\mathbb F)}$ which leaves $GL(n,{\mathbb F})$ stable. It acts on the space of distributions on $GL(n+1,{\mathbb F})$.

Theorem 1 is a Corollary of :
\vskip 0.5cm
{\parindent 0pt
{\bf Theorem $\mathbf 2$:}
{\sl  A distribution on $GL(W)$ which is invariant under the adjoint action of $GL(V)$ is invariant by transposition.}}
\vskip 0.5cm

One can raise a similar question for orthogonal and unitary
groups. Let ${\mathbb D}$ be   either ${\mathbb F}$ or a
quadratic extension of ${\mathbb F}$. If $x\in{\mathbb D}$ then
$\overline x$ is the conjugate of $x$ if ${\mathbb D}\neq{\mathbb
F}$ and is equal to $x$ if ${\mathbb D}={\mathbb F}$.

Let $W$ be a vector space over ${\mathbb D}$ of finite dimension $n+1\geqslant 1$.  Let $\langle  .,.\rangle  $ be a non degenerate hermitian form on $W$. This form is bi-additive and
\[ \langle  dw,d'w'\rangle  =d\,\,\overline{d'}\langle  w,w'\rangle   ,\quad \langle  w',w\rangle  =\overline{\langle  w,w'\rangle  }.\]
Given a ${\mathbb D}-$linear map $u$ from $W$ into itself, its adjoint $u^*$ is defined by the usual formula
\[ \langle  u(w),w'\rangle  =\langle  w,u^*(w')\rangle . \]

Choose a vector $e$ in $W$ such that $\langle  e,e\rangle  \ne 0$; let $U={\mathbb D}e$ and $V=U^{\perp}$ the orthogonal complement. Then $V$ has dimension $n$ and the restriction of the hermitian form to $V$ is non degenerate.

Let $M$ be the unitary group of $W$, that is to say the group of
all ${\mathbb D}-$linear maps $m$ of $W$ into itself which
preserve the hermitian form or equivalently such that $mm^*=1$.
Let $G$ be the unitary group of $V$. With the p-adic topology both
groups are of type lctd (locally compact, totally discontinuous)
and countable at infinity. They are reductive groups of classical
type.

The group $G$ is naturally imbedded into $M$. \vskip 0.5cm

{\parindent 0pt
{\bf Theorem $\mathbf 1'$:}
 {\sl If $\pi $ (resp $\rho $) is an irreducible admissible representation of $M$ (resp of $G$) then
\[ \mathrm{dim}\left  ({\mathrm {Hom}}_G(\pi _{| G},\rho )\right  )\leq 1.\]
}}
Choose a basis $e_1,\dots e_n$ of $V$ such that $\langle  e_i,e_j\rangle  \in{\mathbb F}$. For
\[ w=x_0e+\sum _1^nx_ie_i\]
put
\[ \overline w=\overline x_0e+\sum _1^n\overline {x_i}\, e_i.\]
If $u$ is a ${\mathbb D}-$linear map from $W$ into itself, let $\overline u$ be defined by
\[ \overline u(w)=\overline{u(\overline w)}.\]

Let $\sigma $ be the anti-involution $\sigma (m)=\overline m^{-1}$
of $M$;  Theorem 1' is a consequence of \vskip 0.5cm

{\parindent 0pt
{\bf Theorem $\mathbf 2'$:}
{\sl A distribution on $M$ which is invariant under the adjoint action of $G$ is invariant under $\sigma $.
}}
\vskip 0.5cm

Let us describe briefly our proof. In section 2 we recall why Theorem 2 (2') implies Theorem 1(1').

Then we proceed with $\mathrm{GL}(n)$. The proof is by induction
on  $n$; the case $n=0$ is trivial. In general we first linearize
the problem by replacing the action of $G$ on $\mathrm{GL}(W)$ by
the action on the Lie algebra of $\mathrm{GL}(W)$. As a $G-$module
this Lie algebra is isomorphic to a direct sum $\mathfrak{g}\oplus
V\oplus V^*\oplus\mathbb{F}$ with $\mathfrak{g}$ the Lie algebra
of $G$, $V^*$ the dual space of $V$. The group $G=\mathrm{GL}(V)$
acts trivially on $\mathbb{F}$, by the adjoint action on its Lie
algebra and the natural actions on $V$ and $V^*$. The component
$\mathbb{F}$ plays no role. Let $u$ be a linear bijection of $V$
onto $V^*$ which transforms some basis of $V$ into its dual basis.
The involution may be taken as
$$(X,v,v^*)\mapsto (u^{-1}{\,}^tX\, u,u^{-1}(v^*),u(v)).$$
We have to show that a distribution $T$ on $\mathfrak{g}\oplus V\oplus V^*$ which is invariant under $G$ and skew relative to the involution is 0.

In section 2 we prove that such a distribution must have a
singular support. On the $\mathfrak{g}$ side, using Harish-Chandra
descent we get that the support of $T$ must be contained in
$\mathfrak{z}\times \mathcal{N}\times (V\oplus V^*)$ where
$\mathfrak{z}$ is the center of $\mathfrak{g}$ and $\mathcal{N}$
the cone of nilpotent elements in $\mathfrak{g}$. On the $V\oplus
V^*$ side we show that the support must be contained in
$\mathfrak{g}\times\Gamma$ where $\Gamma$ is the cone $\langle v,
v^*\rangle = 0$ in $V\oplus V^*$. On $\mathfrak{z}$ the action is
trivial so we are reduced to the case of a distribution on
$\mathcal{N}\times\Gamma$.

In section 3 we consider such distributions. The end of the proof
is based on two remarks. First, viewing the distribution as a
distribution on $\mathcal{N}\times (V\oplus V^*)$ its partial
Fourier transform relative to $V\oplus V^*$ has the same
invariance properties and hence must also be supported by
$\mathcal{N}\times \Gamma$. This implies in particular a
homogeneity condition on $V\oplus V^*$. The idea of using Fourier
transform in this kind of situation goes back at least to
Harish-Chandra (\cite{H}) and is conveniently expressed using a
particular case of the Weil or oscillator representation.

For $(v,v^*)\in\Gamma$, let $X_{v,v^*}$ be the map $x\mapsto \langle x,v^*\rangle v$ of $V$ into itself. The second remark is that the one parameter group of transformations
$$(X,v,v^*)\mapsto (X+\lambda X_{v,v^*},v,v^*)$$
is a group of (non linear) homeomorphisms of
$[\mathfrak{g},\mathfrak{g}]\times\Gamma$ which commute with $G$
and the involution. It follows that the image of the support of
our distribution must also be singular. Precisely this allows us
to replace the condition $\langle v,v^*\rangle =0$ by the stricter
condition $X_{v,v^*}\in \mathrm{Im}\,\,\mathrm{ad}X$.

Using the stratification of $\mathcal{N}$ we proceed one nilpotent
orbit at a time, transferring the problem to $V\oplus V^*$ and a
fixed nilpotent matrix $X$. The support condition turns out to be
compatible with direct sum so that it is enough to consider the
case of a principal nilpotent element. In this last situation the
key is the homogeneity condition coupled with an easy induction.

The orthogonal and unitary cases are proved roughly in the same way. In section 4 we reduce the support to the singular set. Here the main difference is that we use Harish-Chandra descent directly on the group. Note that some Levi subgroups have components of type $\mathrm{GL}$ so that theorem 2 has to be assumed. Finally in section 5 we consider the case of a distribution with singular support; the proof follows the same line as in section 3.

We systematically use two classical results : Bernstein's
localization principle and a variant of Frobenius reciprocity
which we call Frobenius descent. For the convenience of the reader
they are both recalled in a short Appendix.

Similar theorems should be true in the archimedean case. A partial
result is given by \cite{AGS}.

Let us add some comments on the Theorems themselves. First note
that Theorem 1' gives an independent proof of a well known theorem
of Bernstein: choose a basis $e_1,\dots ,e_n$ of $V$, add some
vector $e_0$ of $W$ to obtain a basis of $W$ and let $P$ be the
isotropy of $e_0$ in $\mathrm{GL}(W)$. Then Theorem B of
\cite{Ber} says that a distribution on $\mathrm{GL}(W)$ which is
invariant under the action of $P$ is invariant under the action of
$\mathrm{GL}(W)$. Now, by Theorem 1' such a distribution is
invariant under the adjoint action of the transpose of $P$ and the
 group of inner automorphisms is generated by the images of $P$ and its transpose.
We thus get an independent proof of Kirillov conjecture in the
characteristic 0, non archimedean case.

The occurrence of involutions in multiplicity  at most one
problems is of course nothing new. The situation is fairly simple
when all the orbits are stable by the involution thanks to
Bernstein's localization principle and constructibility theorem
(\cite{BZ},\cite{GK}). In our case this is not true : only generic
orbits are stable. Non stable orbits may carry invariant measures
but they do not extend to the ambient space (a similar situation
is already present in \cite{Ber}).

 An illustrative example is the case $n=1$ for $\mathrm{GL}$.
It reduces to $\mathbb{F}^*$ acting on $\mathbb{F}^2$ as
$(x,y)\mapsto (tx,t^{-1}y)$. On the $x$ axis the measure
$d^*x=dx/|x|$ is invariant but does not extend invariantly.
However the symmetric measure
$$f\mapsto \int_{\mathbb{F}^*}f(x,0)d^*x+\int_{\mathbb{F}^*}f(0,y)d^*y$$
does extend.

As in similar cases (for example \cite{JR}) our proof does not
give a simple explanation of why all invariant distributions are
symmetric. The situation would be much better if we had some kind
of density theorem. For example in the $\mathrm{GL}$ case let us
say that an element $(X,v,v^*)$ of $\mathfrak{g}\oplus V\oplus
V^*$ is regular if $(v,Xv,\dots X^{n-1}v)$ is a basis of $V$ and
$(v^*,\dots ,{}^tX^{n-1}v^*)$ is basis of $V^*$. The set of
regular elements is a non empty Zariski open subset; regular
elements have trivial isotropy subgroups. The regular orbits are
the orbits of the regular elements; they are closed, separated by
the invariant polynomials and stable by the involution (see
\cite{RS1}). In particular they carry invariant measures which,
the orbits being closed, do extend and are invariant by the
involution. It is tempting to conjecture that the subspace of the
space of invariant distributions generated by these measures is
weakly dense. This would provide a better understanding of Theorem
2. Unfortunately if true at all, such a density theorem is likely
to be much harder to prove.

Assuming multiplicity at most one, a more difficult question is to find when it is one. Some partial results are known.

For the orthogonal group (in fact the special orthogonal group)
this question has been studied by B. Gross and D.Prasad
(\cite{GP},\cite{P2}) who formulated a precise conjecture. An up
to date account  is given by B.Gross and M.Reeder (\cite{GR}).
In a different setup, in their work on "Shintani" functions
A.Murase and T.Sugano obtained complete results for $GL(n)$ and
the split orthogonal case but only for spherical representations
(\cite{KSM},\cite{MS}). Finally we should mention, Hakim's
publication \cite{Hi}, which, at least for the discrete series,
could perhaps lead to a  different kind of proof.

\subsection*{Acknowledgements}
The first two authors would like to thank their teacher {\bf Joseph
Bernstein} for their mathematical education. They cordially thank
{\bf Joseph Bernstein} and {\bf Eitan Sayag} for guiding them
through this project. They would also like to thank {\bf Yuval
Flicker}, {\bf Erez Lapid}, {\bf Omer Offen} and {\bf Yiannis
Sakellaridis} for useful remarks.

The first two authors worked on this project while participating in
the program \emph{Representation theory, complex analysis and
integral geometry} of the Hausdorff Institute of Mathematics (HIM)
at Bonn joint with Max Planck Institute fur Mathematik. They wish to
thank the organizers of the activity and the director of HIM for
inspiring environment and perfect working conditions.

Finally, the first two authors wish to thank {\bf Vladimir
Berkovich}, {\bf Stephen Gelbart}, {\bf Maria Gorelik} and {\bf
Sergei Yakovenko} from the Weizmann Institute of Science for their
encouragement, guidance and support.

The last author thanks the Math Research Institute of Ohio State
University in Columbus for several invitations which allowed him to
work with the third author.
\section{\hglue -16pt . Theorem 2(2') implies Theorem 1(1')}

A group of type lctd is  a locally compact, totally discontinuous
group which is countable at infinity. We consider  smooth
representations of such groups. If $(\pi ,E_\pi )$ is such a
representation then $(\pi ^*,E_\pi ^*)$ is the smooth
contragradient. Smooth induction is denoted by $Ind$ and compact
induction by $ind$. For any topological space $T$ of type lctd,
$\mathcal {S}(T)$ is the space of functions locally constant,
complex valued, defined on $T$ and with compact support. The space
$\mathcal {S'}(T)$ of distributions on $T$ is the dual space to
$\mathcal {S}(T)$.
\begin{proposition}
 Let $M$ be a lctd group and $N$ a closed subgroup, both unimodular. Suppose that there exists an involutive anti-automorphism $\sigma $ of $M$ such that $\sigma (N)=N$ and such that any distribution on $M$, biinvariant under $N$, is fixed by $\sigma $. Then, for any irreducible admissible representation $\pi $ of $M$
\[ {\mathrm {dim}} \left  ({\mathrm {Hom}} _ M(ind_ N^M(1),\pi )\right  )\times {\mathrm {dim}} \left  ({\mathrm {Hom}} _ M(ind_ N^M(1),\pi ^*)\right  )\leq 1.\]
\end{proposition}
This is well known (see for example \cite{P2}).

{\bf Remark\pointir} There is a variant for the non unimodular case; we will not need it.

\begin{corollary}
 Let  $M$ be a lctd group and $N$ a closed subgroup, both unimodular. Suppose that there exists an involutive anti-automorphism $\sigma $ of $M$ such that $\sigma (N)=N$ and such that any distribution on $M$, invariant under the adjoint action of  $N$, is fixed by $\sigma $. Then, for any irreducible admissible representation $\pi $ of $M$ and any irreducible admissible representation $\rho $ of $N$
\[ \dim\left  ({\mathrm {Hom}}_ N(\pi _{|N},\rho ^*)\right  )\times \dim\left  ({\mathrm {Hom}}_ N((\pi ^*)_{|N},\rho )\right  )\leq 1.\]
\end{corollary}
{\it Proof.} Let $M'=M\times N$ and $N'$ be the closed subgroup of
$M'$ which is the image of the homomorphism $n\mapsto (n,n)$  of
$N$ into $M$. The map $(m,n)\mapsto mn^{-1}$ of $M'$ onto $M$
defines a homeomorphism of $M'/N'$ onto $M$. The inverse map is
$m\mapsto (m,1)N'$. On $M'/N'$ left translations by $N'$
correspond to the adjoint action of $N$ onto $M$. We have a
bijection between the space of distributions $T$ on $M$ invariant
under the adjoint action of $N$ and the space of distributions $S$
on $M'$ which are biinvariant under $N'$. Explicitly
\[ \langle  S,f(m,n)\rangle  =\langle  T,\int _Nf(mn,n)dn\rangle . \]
Suppose that $T$ is invariant under $\sigma $ and consider the involutive anti-automorphism $\sigma '$ of $M'$ given by $\sigma '(m,n)=(\sigma (m),\sigma (n))$. Then
\[ \langle  S,f\circ \sigma '\rangle  =\langle  T,\int _Nf(\sigma (n)\sigma (m),\sigma (n))dn\rangle.  \]
Using the invariance under $\sigma $ and for the adjoint action of $N$ we get
 \begin{eqnarray*}
 \langle  T,\int _Nf(\sigma (n)\sigma (m),\sigma (n))dn\rangle  &=&\langle  T,\int _Nf(\sigma (n)m,\sigma (n))dn\rangle  \\
&=&\langle  T,\int _Nf(mn,n)dn\rangle  \\
&=&\langle  S,f\rangle .
\end{eqnarray*}
\overfullrule=0pt Hence $S$ is invariant under $\sigma '$.
Conversely if $S$ is invariant under $\sigma '$ the same
computation shows that $T$ is invariant under $\sigma $. Under the
assumption of the corollary we can now apply Proposition 1-1 and
we obtain the inequality
 \[ {\mathrm {dim}} \left  ({\mathrm {Hom}} _ {M'}(ind_ {N'}^{M'}(1),\pi \otimes \rho )\right  )\times {\mathrm {dim}} \left  ({\mathrm {Hom}} _ {M'}(ind_ {N'}^{M'}(1),\pi ^*\otimes \rho ^*)\right  )\leq 1.\]
 We know that $Ind_ {N'}^{M'}(1)$ is the smooth contragredient representation of $ind_{N'}^{M'}(1)$; hence
 \[ {\mathrm {Hom}} _ {M'}(ind_ {N'}^{M'}(1),\pi ^*\otimes \rho ^*)\approx {\mathrm{ Hom}} _ {M'}(\pi \otimes \rho , Ind _ {N'}^{M'}(1)).\]
 Frobenius reciprocity tells us that
 \[ {\mathrm {Hom}} _ {M'}\bigl (\pi \otimes \rho , Ind _ {N'}^{M'}(1)\bigr )\approx {\mathrm {Hom}} _ {N'}\bigl ((\pi \otimes \rho )_{|N'},1\bigr ).\]
Clearly
\[ {\mathrm {Hom}} _ {N'}\bigl ((\pi \otimes \rho )_{|N'},1\bigr )\approx {\mathrm {Hom}} _ {N}\bigl (\rho ,(\pi _ {|N})^*\bigr )\approx {\mathrm {Hom}}_ N(\pi _{|N},\rho ^*).\]
Using again Frobenius reciprocity we get
\[ {\mathrm {Hom}} _ {N}\bigl (\rho ,(\pi _ {|N})^*\bigr )\approx {\mathrm {Hom}} _ {M}\bigl (ind_ N^M(\rho ),\pi ^*\bigr ).\]
In the above computations we may replace $\rho $ by $\rho ^*$ and $\pi $ by $\pi ^*$. Finally
\begin{eqnarray*}
{\mathrm {Hom}} _ {M'}(ind_ {N'}^{M'}(1),\pi ^*\otimes \rho ^*)&\approx& {\mathrm {Hom}} _ {N}(\rho ,(\pi _ {|N})^*)\\
&\approx& {\mathrm {Hom}}_ N(\pi _{|N},\rho ^*)\\
&\approx& {\mathrm {Hom}} _ {M}(ind_ N^M(\rho ),\pi ^*).\\
{\mathrm {Hom}} _ {M'}(ind_ {N'}^{M'}(1),\pi \otimes \rho )
&\approx&{\mathrm {Hom}} _ {N}(\rho ^*,((\pi ^*)_ {|N})^*)\\
&\approx& {\mathrm {Hom}}_ N((\pi ^*)_{|N},\rho )\\
&\approx& {\mathrm {Hom}} _ {M}(ind_ N^M(\rho ^*),\pi ).
\end{eqnarray*}\hfill$\Box$\smallskip

Going back to our situation and keeping the notations of the
introduction consider first the case of the general linear group.
We take $M=GL(W)$ and $N=GL(V)$. Let $E_\pi $ be the space of the
representation $\pi $  and let $E_\pi ^*$ be the smooth  dual
(relative to the action of $GL(W))$. Let $E_\rho $ be the space of
$\rho $ and $E^*_\rho $ be the smooth dual for the action of
$GL(V)$. We know (\cite{BZ} section 7), that the contragredient
representation $\pi ^*$ in $E_\pi ^*$ is isomorphic to the
representation $g\mapsto \pi (^tg^{{-1}})$ in $E_\pi $. The same
is true for $\rho ^*$. Therefore an element of ${\mathrm {Hom}}_
N(\pi _{|N},\rho ^*)$ may be described as a linear map $A$ from
$E_\pi $ into $E_\rho $ such that, for $g\in N$
\[ A\pi (g)=\rho (^tg^{-1})A.\]
An element of ${\mathrm {Hom}}_ N((\pi ^*)_{|N},\rho )$ may be described as a linear map $A'$ from $E_\pi $ into $E_\rho $ such that, for $g\in N$
\[ A'\pi (^tg^{-1})=\rho (g)A'.\]
We have obtained the same set of linear maps:
\[ {\mathrm {Hom}}_ N((\pi ^*)_{|N},\rho )\approx {\mathrm {Hom}}_ N(\pi _{|N},\rho ^*).\]
We are left with 2 possibilities: either both spaces have dimension 0 or they both have dimension 1 which is exactly what we want.

From now on we forget Theorem 1 and  prove Theorem  2.

Consider the orthogonal/unitary case, with the notations of the introduction.
In Chapter 4 of \cite{MVW} the following result is proved. Choose $\delta \in GL_{\mathbb F}(W)$ such that $\langle  \delta w,\delta w'\rangle  =\langle  w',w\rangle  $. If $\pi $ is an irreducible admissible representation of $M$, let $\pi ^*$ be its smooth contragredient and define $\pi ^\delta $ by
\[ \pi ^\delta (x)=\pi (\delta x\delta ^{-1}).\]
Then $\pi ^\delta $ and $\pi ^*$ are equivalent. We choose $\delta =1$ in the orthogonal case ${\mathbb D}={\mathbb F}$. In the unitary case, fix an orthogonal basis of $W$, say $e_1,\dots ,e_{n+1}$, such that $e_2,\dots ,e_{n+1}$ is a basis of $V$; put $\langle  e_i,e_i\rangle  =a_i$. Then
\[ \langle  \sum x_ie_i,\sum y_je_j\rangle  =\sum a_ix_i\overline{y_i}.\]
Define $\delta $ by
\[ \delta \left  (\sum x_ie_i\right  )=\sum \overline{x_i}e_i.\]
Note that $\delta ^2=1$.

Let $E_\pi $ be the space of $\pi $. Then, up to equivalence, $\pi ^*$ is the representation $m\mapsto \pi (\delta m\delta ^{-1})$. If $\rho $ is an admissible irreducible representation of $G$ in a vector space $E_\rho $ then an element $A$ of ${\mathrm {Hom}}\left  (\pi ^*_{|G},\rho \right  )$ is a linear map from $E_\pi $ into $E_\rho $ such that
\[ A\pi (\delta g\delta ^{-1})=\pi (g)A,\quad g\in G.\]
In turn the contragredient $\rho ^*$ of $\rho $  is equivalent to the representation $g\mapsto \rho (\delta g\delta ^{-1})$ in $E_\rho $. Then an element $B$ of ${\mathrm {Hom}}\left  (\pi _{|G},\rho ^*\right  )$ is a linear map from $E_\pi $ into $E_\rho $ such that
\[ B\pi (g)=\rho (\delta g\delta ^{-1})B,\quad g\in G.\]
As $\delta ^2=1$ the conditions on $A$ and $B$ are the same:
\[ {\mathrm {Hom}}\left  (\pi ^*_{|G},\rho \right  )\approx {\mathrm {Hom}}\left  (\pi _{|G},\rho ^*\right  ). \]
However, assuming Theorem 2', by Corollary 1-1 we have
\[ \dim\biggl ({\mathrm {Hom}}\left  (\pi ^*_{|G},\rho \right  )\biggr )\times \dim\biggl ({\mathrm {Hom}}\left  (\pi _{|G},\rho ^*\right  )\biggr )\leq 1.
\]
so that both dimensions are 0 or 1. Replacing $\rho $ by $\rho ^*$ we get Theorem 1'. From now on we forget about Theorem 1'.

\section{\hglue -16pt . Reduction to the singular set : the  GL(n) case}

If $ H$  is a topological group of type lctd, acting continuously
on a topological space $ E$  of the same type and if $\chi $  is a
continuous character of $ H$  we denote by $ {\mathcal
S}'(E)^{H,\chi }$  the space of distributions $ T$  on $ E$  such
that $ \langle T,f(h^{-1}x)\rangle =\chi (h)\langle T,f\rangle $  for
any $ f\in  {\mathcal  S}(E)$  and any $ h\in  H$.

Consider the case of the general linear group. From the decomposition $ W=V\oplus {\mathbb   F}e$  we get, with obvious identifications
$$ {\mathrm {End}}  (W)={\mathrm {End}}  (V)\oplus V\oplus V^*\oplus {\mathbb   F}.$$
Note that $ {\mathrm {End}}  (V)$  is the Lie algebra $ {\mathfrak
g}$  of $ G$. The group $ G$  acts on $ {\mathrm {End}}  (W)$  by
$ g(X,v,v^*,t)=(gXg^{-1},gv,{}^t g^{-1}v^*,t)$. As before choose a
basis $ (e_1,\dots ,e_n)$  of $ V$  and let $ (e_1^*,\dots
,e_n^*)$  be the dual basis of $ V^*$. Define an isomorphism $ u$
of $ V$  onto $ V^*$  by $ u(e_i)=e_i^*$. On $ GL(W)$  the
involution $ \sigma $  is $ h\mapsto u^{-1}{}^th^{-1}u$. It
depends upon the choice of the basis but the action on the space
of invariant distributions does not depend upon this choice.

 It will be convenient to introduce an extension  $ \widetilde  G$  of $ G$. Let $ {\mathrm {Iso}}(V,V^*)$  be the set of isomorphisms of $ V$  onto $ V^*$. We define $ \widetilde  G=G\cup  \mathrm {Iso}(V,V^*)$. The group law, for $ g,g'\in  G$  and $ u,u'\in   \mathrm {Iso}(V,V^*)$  is
$$ g \times g'=gg',\,\,\, u \times g=ug,\,\,\, g \times u=^t\hskip -3 pt g^{-1}u,\,\,\, u \times u'=^t\hskip -3 pt u^{-1}u'.$$
Now from $ W=V\oplus {\mathbb   F}e$  we obtain an identification
of the dual space $ W^*$  with $ V^*\oplus {\mathbb   F}e^*$  with
$ \langle e^*,V\rangle =(0)$  and $ \langle e^*,e\rangle =1$. Any
$ u$  as above extends to an isomorphism of $ W$  onto $ W^*$  by
defining $ u(e)=e^*$. The group $ \widetilde  G$  acts on $ GL(W)$
:
$$ h\mapsto ghg^{-1},\,\,\,\, h\mapsto ^t\hskip -3pt (uhu^{-1}),\quad g\in  G,\, h\in  GL(W),\,\, u\in   \mathrm {Iso}(V,V^*)$$
and also on $ {\mathrm {End}}  (W)$  with the same formulas.

Let $\chi $  be the character of $ \widetilde  G$  which is 1 on $
G$  and $ -1$  on $  \mathrm {Iso}(V,V^*)$. Our goal is to prove
that $ {\mathcal  S}'(GL(W))^{\widetilde  G, \chi  }=(0)$.
\begin{proposition}
If $ {\mathcal  S}'({\mathfrak  g}\oplus V\oplus V^*)^{\widetilde
G, \chi  }=(0)$  then $ {\mathcal  S}'(GL(W))^{\widetilde  G, \chi
}=(0)$.
\end{proposition}
 {\it Proof.} We have $ {\mathrm {End}}(W)=\bigl (\,\,{\mathrm {End}}  (V)\oplus V\oplus V^*\bigr )\oplus {\mathbb   F}$  and the action of $ \widetilde  G$  on $ {\mathbb   F}$  is trivial thus  $ {\mathcal  S}'({\mathfrak  g}\oplus V\oplus V^*)^{\widetilde  G, \chi  }=(0)$  implies that $ {\mathcal  S}'({\mathrm {End}}   (W))^{\widetilde  G, \chi  }=(0).$
Let $ T\in  {\mathcal  S}'(GL(W))^{\widetilde  G, \chi  }$. Let $
h\in  GL(W)$  and choose a compact open neighborhood $ K$  of $
{\mathrm {Det}}\,\,  h$  such that $ 0\notin K$. For $ x\in
{\mathrm {End}  } (W)$  define  $ \varphi (x)=1$  if $  {\mathrm
{Det}}  x\in  K $  and $ \varphi (x)=0$  otherwise. Then $ \varphi
$  is a locally constant function. The distribution $ (\varphi
_{|GL(W)})T$  has a support which is closed in $  {\mathrm {End}}
(W)$  hence may be viewed as a distribution on $  {\mathrm {End} }
(W)$. This distribution belongs to $ {\mathcal  S}'({\mathrm
{End}}   (W))^{\widetilde  G, \chi  }$  so it must be equal to 0.
It follows that $ T$  is 0 in the neighborhood of $ h$. As $ h$
is arbitrary we conclude that $ T=0$.

\hfill$\Box$\smallskip

Our task is now to prove that $ {\mathcal  S}'({\mathfrak g}\oplus
V\oplus V^*)^{\widetilde  G, \chi  }=(0)$. We shall use induction
on the dimension $ n$  of $ V$.  The action of $ \widetilde  G$
is, for $ X\in  {\mathfrak  g},\, v\in  V,\, v^*\in  V^*,\, g\in
G,\,  u\in   \mathrm {Iso}(V,V^*)$
$$ (X,v,v^*)\mapsto (gXg^{-1},gv,^t\hskip -3pt g^{-1}v^*),\,\,\, (X,v,v^*)\mapsto (^t\hskip -1pt (uXu^{-1}),^t\hskip -3pt u^{-1}v^*,uv).$$
The case $n=0$  is trivial.

We suppose that  $ V$   is of dimension $ n\geq 1$, assuming the
result up to dimension $ n-1$ and for all $\mathbb F$. If $ T\in
{\mathcal  S}'({\mathfrak  g}\oplus V\oplus V^*)^{\widetilde  G,
\chi  }$  we are going to show that its support   is contained in
the "singular set". This will be done in two stages.

On $ V\oplus V^*$  let $ \Gamma $  be the cone $ \langle
v^*,v\rangle =0$. It is stable under $ \widetilde  G$.
\begin{lemma}
 The support of $ T$  is contained in $ {\mathfrak  g} \times \Gamma $.
\end{lemma}

{\it Proof.} For $ (X,v,v^*)\in  {\mathfrak  g}\oplus V\oplus V^*$
put $ q(X,v,v^*)=\langle v^*,v\rangle $. Let $ \Omega $  be the
open subset $ q\ne 0$. We have to show that $ {\mathcal S}'(\Omega
)^{\widetilde  G, \chi  }=(0)$. By Bernstein's localization
principle (Corollary 6-1) it is enough to prove that, for any fiber $ \Omega
_t=q^{-1}(t),\,\,\, t\ne 0$, one has $ {\mathcal S}'(\Omega
_t)^{\widetilde  G, \chi  }=(0)$.

 $ G$  acts transitively on the quadric $ \langle v^*,v\rangle =t$. Fix a decomposition
 $ V={\mathbb   F}\varepsilon \oplus V_1$  and identify $ V^*={\mathbb   F}\varepsilon ^*\oplus V_1^*$
 with $ \langle \varepsilon ^*,\varepsilon \rangle =1$. Then $ (X,\varepsilon ,t\varepsilon ^*)\in  \Omega _t$
 and the isotropy subgroup of $ (\varepsilon ,t\varepsilon ^*)$  in $ \widetilde  G$  is, with an obvious notation
 $ \widetilde  G_{n-1}$. By Frobenius descent (Theorem 6-2) there is a linear bijection between
 $ {\mathcal  S}'(\Omega _t)^{\widetilde  G, \chi  }$  and the space
 $ {\mathcal  S}'({\mathfrak  g})^{\widetilde  G_{n-1}, \chi  }$  and this last space is $ (0)$  by induction.

 \hfill$\Box$\smallskip

Let $ {\mathfrak  z}$  be the center of $ {\mathfrak  g}$  that is
to say the space of scalar matrices. Let $ {\mathcal  N}\subset
[{\mathfrak  g},{\mathfrak  g}]$  be the nilpotent cone in $
{\mathfrak  g}$.
\begin{lemma}
 The support of $ T$  is contained in $ {\mathfrak  z} \times {\mathcal  N} \times \Gamma $.
\end{lemma}
{\it Proof.} We use Harish-Chandra's descent.  For $ X\in
{\mathfrak  g}$  let $ X=X_s+X_n$  be the Jordan decomposition of
$ X$  with $ X_s$  semisimple and $ X_n$  nilpotent. This
decomposition commutes with the action of $ \widetilde  G$. The
centralizer $ Z_G(X)$  of an element $ X\in  {\mathfrak  g}$   is
unimodular (\cite{SS} page 235) and there exists an isomorphism $
u$  of $ V$  onto $ V^*$  such that $ ^tX=uXu^{-1}$  (any matrix
is conjugate to its transpose). It follows that the centralizer $
Z_{\widetilde  G}(X)$  of $ X$  in $ \widetilde  G$, a semi direct
product of $ Z_G(X)$  and $ S_2$  is also unimodular.

Let $ E$  be the vector space of monic polynomials, of degree $ n$
, with coefficients in $ {\mathbb   F}$. For $ p\in  E$, let $
{\mathfrak  g}_p$  be the set of all $ X\in  {\mathfrak  g}$  with
characteristic polynomial $ p$. Note that $ {\mathfrak  g}_p$  is
fixed by $ \widetilde  G$. By Bernstein localization principle (Corollary 6-1) it
is enough to prove that if $ p$  is not $ (T-\lambda )^n$  for
some $ \lambda $  then $ {\mathcal  S}'({\mathfrak  g}_p \times V
\times V^*)^{\widetilde  G, \chi  }=(0)$.

Fix $ p$. We claim that the map $ X\mapsto X_s$  restricted to $
{\mathfrak  g}_p$  is continuous. Indeed let $ \widetilde {\mathbb
F}$  be a finite Galois extension of $ {\mathbb   F}$ containing
all the roots of $ p$. Let
$$ p(\xi )=\prod _1^s(\xi -\l_i)^{n_i}$$
be the decomposition of $ p$. Recall that if $ X\in  {\mathfrak
g}_p$  and $ V_i=\mathrm{Ker }   (X-\lambda _I)^{n_i}$  then $
V=\oplus V_i$  and the restriction of $ X_s$  to $ V_i$  is the
multiplication by $ \lambda _i$. Then choose a polynomial $ R$,
with coefficients in $ \widetilde  {\mathbb   F}$  such that for
all $ i$, $ R$  is congruent to $ \lambda _i$  modulo $ (\xi
-\lambda _i)^{n_i}$  and $ R(0)=0$  to respect the tradition.
Clearly $ X_s=R(X)$. As the Galois group of  $ \widetilde {\mathbb
F}$   over $ {\mathbb   F}$  permutes the $ \lambda _i$ we may
even choose $ R\in  {\mathbb   F}[\xi ]$. This implies the
required continuity.

There is only one semi-simple orbit $ \gamma  _p$  in $ {\mathfrak
g}_p$  and it is closed. We use Frobenius descent  for the map $
(X,v,v^*)\mapsto X_s$  from $ {\mathfrak  g}_p \times V \times
V^*$  to $ \gamma  _p$.

Fix $ a\in  \gamma  _p$ ;  its fiber is the product of $ V\oplus
V^*$  by the set of nilpotent elements which commute with $ a$. It
is a closed subset of the centralizer  $ {\mathfrak m}={\mathfrak
Z}_{{\mathfrak  g}}(a)$  of $ a$  in $ {\mathfrak g}$. Let $
M=Z_G(a)$  and $ \widetilde  M=Z_{\widetilde  G}(a)$.

Following (\cite{SS} ) let us describe these centralizers.
 Let $P$ be the minimal polynomial of $ a$ ; all its roots are simple. Let $P=P_1\dots P_r$ be
the decomposition of $P$ into irreducible factors, over ${\mathbb    F}$. Then the $P_i$
 are two by two relatively prime. If $V_i={\mathrm Ker}    P_i(a)$, then $V=\oplus   V_i$ and $V^*=\oplus   V_i^*$. An element $x$ of $G$
which commutes with $a$ is given by  a family $\{x_1,\dots ,x_r\}$
where each $x_i$ is a linear map from $V_i$ to $V_i$, commuting
with the restriction of $a$ to $V_i$. Now ${\mathbb    F}[\xi ]$
acts on $V_i$, by specializing $\xi $ to $a_{|V_i}$ and $P_i$ acts
trivially so that, if ${\mathbb    F}_i={\mathbb    F}[\xi
]/(P_i)$, then $V_i$ becomes a vector space over ${\mathbb
F}_i$. The ${\mathbb    F}-$linear map $x_i$ commutes with $a$ if
and only if it is ${\mathbb F}_i-$linear.\hfill\break

Fix $i$. Let $\ell$ be a non zero ${\mathbb    F}-$linear form on
${\mathbb    F}_i$. If $v_i\in V_i$ and $v'_i\in V_i^*$ then
$\lambda \mapsto \langle  \lambda v_i,v'_i\rangle  $ is an
${\mathbb    F}-$linear form on ${\mathbb    F}_i$, hence there
exists a unique element $S(v_i,v'_i)$ of ${\mathbb    F}_i$ such
that $\langle  \lambda v_i,v'_i\rangle  =\ell\left  (\lambda
S(v_i,v'_i)\right  )$. One checks trivially that $S$ is ${\mathbb
F}_i-$linear with respect to each variable and defines a non
degenerate duality, over ${\mathbb    F}_i$ between $V_i$ and
$V_i^*$. Here ${\mathbb    F}_i$ acts on $V_i^*$ by transposition,
relative to the ${\mathbb    F}-$duality $\langle  .,.\rangle $,
of the action on $V_i$. Finally if $x_i\in  \mathrm {End}
_{{\mathbb    F}_i} V_i$, its transpose, relative to the duality
$S(.,.)$ is the same as its transpose relative to the duality
$\langle   .,.\rangle  $.

Thus $M$ is a product of linear groups and the situation $(M,V,V^*)$ is a composite case, each component being
a linear case (over various extensions of ${\mathbb    F}$).

Let $ u$  be an isomorphism of $ V$  onto $ V^*$  such that $
^ta=uau^{-1}$  and that, for each $ i,\,\,\, u(V_i)=V_i^*$. Then $
u\in  \widetilde  M$  and $ \widetilde  M=M\cup uM$.

Suppose that $ a$  does not belong to the center of $ {\mathfrak
g}$. Then each $ V_i$  has dimension strictly smaller than $ n$
and we can use the inductive assumption. Therefore $ {\mathcal
S}'({\mathfrak  m}\oplus V\oplus V^*)^{\widetilde  M, \chi
}=(0)$. However the nilpotent cone $ {\mathcal  N}_{\mathfrak  m}$
in $ {\mathfrak  m}$  is a closed subset so $ {\mathcal
S}'({\mathcal N}_{\mathfrak  m} \times V \times V^*)^{\widetilde
M, \chi }=(0)$  which is what we need.

\hfill$\Box$\smallskip

If $ a$  belongs to the center then $ \widetilde  M=\widetilde  G$
and the fiber is $ (a+{\mathcal  N}) \times V \times V^*$.
Therefore we have proved the following Proposition:

\begin{proposition}
If $ T\in  {\mathcal  S}'({\mathfrak  g}\oplus V\oplus V^*)^{\widetilde  G, \chi  }$ then the support of $T$ is contained in ${\mathfrak z}\times{\mathcal N}\times{\Gamma}$.\hfill\break
 If  $ {\mathcal  S}'({\mathcal  N} \times\Gamma)^{\widetilde  G, \chi  }=(0)$  then  $ {\mathcal  S}'({\mathfrak  g}\oplus V\oplus V^*)^{\widetilde  G, \chi  }=(0)$.
\end{proposition}

{\textbf{Remark}\pointir}Strictly speaking the singular set is
defined as the set of all $(X,v,v^*)$ such that for  any
polynomial $P$ invariant under $\widetilde G$ one has
$P(X,v,v^*)=P(0)$. So we should take care of the invariants
$P(X,v,v^*)=\langle v^*,X^pv\rangle$ for all $p$ and not only for
$p=0$. It can be proved, a priori, that the support of the
distribution $T$ has to satisfy these extra conditions. As this is
not needed in the sequel we omit the proof.

\section{\hglue -16pt . End of the proof for GL(n)}
In this section we consider a distribution $ T\in  {\mathcal
S}'({\mathcal  N} \times \Gamma )^{\widetilde  G, \chi  }$  and
prove that $ T=0$. The following observation will play a crucial
role.

Choose a non trivial additive character $\psi$ of $\mathbb F$. On $V\oplus V^*$ we have the bilinear form
$$\bigl ((v_1,v_1^*),(v_2,v_2^*)\bigr )\mapsto \langle v_1^*,v_2\rangle +\langle v_2^*,v_1\rangle.$$
Define the Fourier transform by
$$\widehat \varphi (v_2,v_2^*)=\int _{V\oplus V^*}\varphi (v_1,v_1^*)\, \psi (\langle v_1^*,v_2\rangle +\langle v_2^*,v_1\rangle )\, dv_1dv_1^*$$
with $dv_1dv_1^*$ is normalized so that there is no constant factor appearing in the inversion formula.

This Fourier transform commutes with the action of $\widetilde G$;
hence the (partial) Fourier transform $\widehat T$ of our
distribution $T$ has the same invariance properties and the same
support conditions as $T$ itself.

Let $ {\mathcal  N}_i$  be the union of nilpotent orbits of
dimension at most $ i$. We will prove, by descending  induction on
$ i$, that the support of any $(\widetilde{G},\chi)-$equivariant
distribution  must be contained in $ {\mathcal  N}_i \times \Gamma
$. Suppose we already know that, for some $ i$, the support must
be contained in $ {\mathcal N}_i \times \Gamma $. We must show
that, for any nilpotent orbit $ {\mathcal  O}$  of dimension $ i$,
the restriction of the distribution to $ {\mathcal  O} \times
\Gamma $ is 0.

If $ v\in  V$  and $ v^*\in  V^*$  we call $ X_{v,v^*}$  the rank
one map $ x\mapsto \langle v^*,x\rangle v$. Let
$$ \nu_\lambda (X,v,v^*)=(X+\lambda X_{v,v^*},v,v^*),\quad (X,v,v^*)\in  {\mathfrak  g} \times \Gamma ,\,\,\, \lambda \in  {\mathbb   F}.$$
Then $ \nu_\lambda $  is a one parameter group of homeomorphisms
of $ {\mathfrak  g} \times \Gamma $  and note that $ [{\mathfrak
g},{\mathfrak  g}] \times \Gamma $  is invariant. The key
observation is that  {\it $ \nu_\lambda $  commutes with the
action of $ \widetilde  G$ }. Therefore the image of $ T$  by $
\nu_\lambda $  transforms according to the character $\chi $  of $
\widetilde  G$. Its support is contained in $ [{\mathfrak
g},{\mathfrak  g}] \times \Gamma $ and hence must be contained in $
{\mathcal  N} \times \Gamma $  and in fact in $ {\mathcal  N}_i
\times \Gamma $. This means that if $ (X,v,v^*)$  belongs to the
support of $ T$  then, for all $ \lambda $, $ (X+\lambda
X_{v,v^*},v,v^*)$  must belong to $ {\mathcal  N}_i \times \Gamma
$.

The orbit $ {\mathcal  O}$  is open in $ {\mathcal  N}_i$.  Thus if
$ X\in  {\mathcal  O}$  the condition $ X+\lambda X_{v,v^*}\in
{\mathcal  N}_i$  implies that, at least for $ |\lambda |$  small
enough, $ X+\lambda X_{v,v^*}\in  {\mathcal  O}$. It follows that
$ X_{v,v^*}$  belongs to the tangent space to $ {\mathcal  O}$  at
the point $ X$  and this tangent space is the image of $
\mathrm{ad }   X$.

Let us call $Q(X)$ the set of all pairs $(v,v^*)$ such $X_{v,v^*}\in \mathrm {Im\,\, ad}X$.

 Therefore it is enough to prove the following Lemma:
\begin{lemma}
 Let $ T\in  {\mathcal  S}'({\mathcal  O} \times V \times V^*)^{\widetilde  G, \chi  }$. Suppose that the support of $ T$  and of $ \widehat T$  are contained    in the set of triplets  $ (X,v,v^*)$  such that $ (v,v^*)\in  Q(X)$.
Then $ T=0$.
\end{lemma}
Note that the trace of $ X_{v,v^*}$  is $ \langle v^*,v\rangle $
and that $ X_{v,v^*}\in  {\mathrm Im}\,\, \mathrm{ad } \,  X$
implies that its trace is 0. Therefore $ Q(X)$  is contained in $
\Gamma $.

We proceed in three steps. First we transfer the problem to $
V\oplus V^*$  and a fixed nilpotent endomorphism $ X$. Then we
show that if  Lemma 3-1 is true for $ (V_1,X_1)$  and $ (V_2,X_2)$
then it is true for the direct sum $ (V_1\oplus V_2,X_1\oplus
X_2)$. Finally using the decomposition of $ X$  in Jordan blocks
we are left with the case of a principal nilpotent element for
which we give a direct proof, using Weil representation.

Consider the map $ (X,v,v^*)\mapsto X$  from $ {\mathcal  O}
\times V \times V^*$  onto $ {\mathcal  O}$. Choose $X\in{\mathcal
O}$ and let $ C$ (resp $ \widetilde  C$)  be the stabilizer in $
G$ (resp. in $ \widetilde  G$)  of an element $ X$  of $ {\mathcal
O}$ ; both groups are unimodular, hence we may use Frobenius
descent (Theorem 6-2).

 Now we have to deal with a distribution, which we still call $T$, which belongs to $   {\mathcal  S}'(V\oplus V^*)^{\widetilde  C, \chi  }$  such that both $ T$  and its Fourier transform are supported by $Q(X)$. Let us say that $ X$  is nice if the only such  distribution is 0. We want to prove that all nilpotent endomorphisms are nice.

 \begin{lemma}
  Suppose that we have a decomposition $ V=V_1\oplus V_2$  such that $ X(V_i)\subset V_i$. Let $ X_i$  be the restriction of $ X$  to $ V_i$. Then if $ X_1$  and $ X_2$  are nice, so is $ X$.
 \end{lemma}
 {\it Proof of Lemma 3.2.} Let $ Q(X)$  be the set of pairs $ (v,v^*)$  such that $ X_{v,v^*}$  belongs to the image of $ \mathrm{ad }   X$. Let $ (v,v^*)\in  Q(X)$  and choose $ A\in  {\mathfrak  g}$  such that $ X_{v,v^*}=[A,X]$. Decompose $ v=v_1+v_2,\,\, v^*=v_1^*+v_2^*$  and put
 $$ A=\pmatrix{A_{1,1}&A_{1,2}\cr A_{2,1}&A_{2,2}\cr}.$$
 Writing $ X_{v,v^*}$  as a 2 by 2 matrix and looking at the diagonal blocks one gets that $ X_{v_i,v_i^*}=[A_{i,i},X_i]$. This means that
 $$ Q(X)\subset Q(X_1) \times Q(X_2).$$
 For $ i=1,2$  let $ C_i$  be the centralizer of $ X_i$  in $ GL(V_i)$  and $ \widetilde  C_i$  the corresponding extension by $ S_2$.
 Let $ T$  be a distribution as above and let $ \varphi _2\in  {\mathcal  S}(V_2\oplus V_2^*)$. Let $ T_{1}$  be the distribution on $ V_1\oplus V_1^*$  defined by $ \varphi _1\mapsto \langle T,\varphi _1\otimes \varphi _2\rangle $. The support of $ T_{1}$  is contained in $ Q(X_1)$  and $ T_{1}$  is invariant under the action of  $ C_1$. We have
 $$ \langle \widehat T_{1},\varphi _1\rangle =\langle T_{1},\widehat \varphi _1\rangle =\langle T, \widehat {\varphi _1}\otimes \varphi _2\rangle =\langle \widehat T,\check {\varphi _1}\otimes \widehat{\varphi _2}\rangle .$$
 Here $ \check \varphi _1(v_1,v_1^*)=\varphi _1(-v_1,-v_1^*)$. By assumption the support of $ \widehat T$  is contained in $ Q(X)$  so that the support of $ \widehat{T_1}$  is supported in $ -Q(X_1)=Q(X_1)$. Because $ (X_1)$  is nice this implies that $ T_1$  in invariant under $ \widetilde  C_1$. Imbedding $ \widetilde  C_1$  into $ \widetilde  C$  we get that $ T$  is invariant under $ \widetilde  C_1$.   Similarly it is invariant under $ \widetilde  C_2$. However the subgroup $ \widetilde  C_1 \times \widetilde  C_2$  of $ \widetilde  C$  is not contained in $ C$  so that $ T$  must be invariant under $ \widetilde  C$ and  hence must be 0.

 \hfill$\Box$\hfill

 Decomposing $ X$  into Jordan blocks we still have to prove Lemma 3-1 for a principal nilpotent element. We need some preliminary results.

 \begin{lemma}
 The distribution $ T$  satisfies the following homogeneity condition:
$$ \langle T,f(tv,tv^*)\rangle =|t|^{-n}\langle T,f(v,v^*)\rangle .$$
\end{lemma}
    {\it Proof of Lemma 3.3.} We use a particular case of Weil or oscillator representation. Let $ E$  be a vector space over $ {\mathbb   F}$  of finite dimension $ m$. To simplify assume that $ m$  is even. Let $ q$  be a non degenerate quadratic form on $ E$  and let $ b$  be the bilinear form
    $$ b(e,e')=q(e+e')-q(e)-q(e').$$
Fix a continuous non trivial additive character $ \psi  $  of $
{\mathbb   F}$. We define the Fourier transform on $ E$  by
$$ \widehat f(e')=\int  _Ef(e)\psi  (b(e,e'))de$$
where $ de$  is the self dual Haar measure.

There exists (\cite{RS2}) a representation $ \pi $  of $
SL(2,{\mathbb   F})$  in $ {\mathcal  S}(E)$  such that:
 \begin{eqnarray*}
\pi \pmatrix{1&u\cr 0&1\cr}f(e)&=&\psi  (uq(e))f(e)\\
\pi \pmatrix{t&0\cr 0&t^{-1}}f(e)&=&\frac{\gamma  (q)}{\gamma  (tq)}|t|^{m/2}f(te)\\
\pi \pmatrix{0&1\cr -1&0\cr}f(e)&=&\gamma  (q)\widehat f(e).\\
\end{eqnarray*}
The $ \gamma  (tq)$  are complex numbers of modulus 1. In
particular if $(E,q)$  is a sum of hyperbolic planes these numbers
are all equal to 1.

We have a contragredient action in the dual space $ {\mathcal
S'(E)}$.

Suppose that $ T$  is a distribution on $ E$  such that $ T$  and
$ \widehat T$  are supported by the isotropic cone $ q(e)=0$. This
means that
$$ \langle T,\pi \pmatrix{1&u\cr 0&1\cr}f\rangle =\langle T,f\rangle ,\quad \langle \widehat T,\pi \pmatrix{1&u\cr 0&1\cr}f\rangle =\langle \widehat T,f\rangle . $$
Using the relation
$$ \langle \widehat T,\varphi \rangle =\langle T,\overline {\gamma  (q)}\pi \pmatrix{0&1\cr -1&0\cr}f\rangle $$
the second relation is equivalent to
$$ \langle T,\pi \pmatrix{1&0\cr -u&0\cr}f\rangle =\langle T,f\rangle .$$
The matrices
$$ \pmatrix{1&u\cr 0&1\cr},\quad {\mathrm {and}}\quad \pmatrix{1&0\cr u&1\cr}\quad u\in  {\mathbb   F}$$
generate  the group $ SL(2,{\mathbb   F})$. Therefore the
distribution $ T$  is invariant by $ SL(2,{\mathbb   F)}$. In
particular
$$ \langle T,f(te)\rangle =\frac{\gamma  (tq)}{\gamma  (q)}|t|^{-m/2}\langle T,f\rangle $$
and $ T=\gamma  (q)\widehat T$.

{\bf Remark\pointir } For $ m$   even $ \gamma  (tq)/\gamma  (q)$
is a character and there do exist non zero distributions invariant
under $ SL(2,{\mathbb   F})$. In odd dimension we get a
representation of the 2-fold covering of $ SL(2,{\mathbb   F})$
and we obtain the same homogeneity condition. However $ \gamma
(tq)/\gamma  (q)$  is not a character; hence the distribution $
T$  must be 0.

In our situation we take $ E=V\oplus V^*$  and $ q(v,v^*)=\langle
v^*,v\rangle $. Then
$$ b\Bigl ((v_1,v_1^*),(v_2,v_2^*)\Bigr )=\langle v_1^*,v_2\rangle +\langle v_2^*,v_1\rangle.$$
The Fourier transform commutes with the action of $ \widetilde  G$. Both $
T$  and $ \widehat T$   are supported by $ Q(X)$  which is
contained in $ \Gamma $. As $ \gamma  (tq)=1$  this proves the
Lemma and also that $ T=\widehat T$. \hfill$\Box$\hfill

{\bf Remark\pointir} The same type of argument could have been
used for the quadratic form $ Tr(XY)$  on $ {\mathfrak
s}{\mathfrak  l}(V)=[{\mathfrak  g},{\mathfrak  g}]$. This would
have given a short proof for even $ n$  and a homogeneity
condition for odd $ n$.

Now we find $ Q(X)$.
\begin{lemma}
 If $ X$  is principal then $ Q(X)$  is the set of pairs $ (v,v^*)$  such that for $ 0\leq k<n$, $ \langle v^*,X^kv\rangle =0$.
\end{lemma}
{\it Proof of Lemma 3.4.} Choose a basis $ (e_1,\dots ,e_n)$  of $
V$  such that $ Xe_1=0$  and $ Xe_j=e_{j-1}$  for $ j\geq 2$.
Consider the map $ A\mapsto XA-AX$  from the space of n by n
matrices into itself. A simple computation shows that the kernel
of this map, that is to say the Lie algebra $ {\mathfrak  c}$  of
the centralizer  $ C$,  is the space of polynomials (of degree at
most $ n-1$ ) in $ X$. It is of dimension $ n$. The image is of
codimension $ n$  and calling $ b_{i,j}$  the coefficients of an n
by n matrix, a set of independent equations for this image  is
$$ \sum _{j=1}^{n-r}b_{j+r,j}=0,\quad r=0,\dots n-1.$$
Let $ (e_1^*,\dots ,e_n^*)$  be the dual basis. Call $ x_1,\dots
,x_n$  the coordinates of $ v$  and $ (x_1^*,\dots ,x_n^*)$  the
coordinates of $ v^*$. The matrix of $ X_{v,v^*}$  is then given
by $ b_{i,j}=x_ix_j^*$  and we get the lemma.

\hfill$\Box$\hfill

{\it End of the proof of Lemma 3.1. }For $X$ principal, we proceed
by induction on $ n$. Keep the above notations. The centralizer $
C$  of $ X$  is the space of polynomials (of degree at most $ n-1$
) in $ X$  with non zero constant term. In particular the orbit $
\Omega $  of $ e_n$  is the open subset $ x_n\ne 0$. We shall
prove that the restriction of $ T$  to $ \Omega  \times V^*$ is 0.
Note that the centralizer of $ e_n$  in $ C$  is trivial. By
Frobenius descent (Theorem 6-2), to the restriction of $ T$  corresponds a
distribution $ R$  on $ V^*$  with support in the set of $ v^*$
such that $ (e_n,v^*)\in  Q(X)$. By the last Lemma this means that
$ R$  is a multiple $ a\delta $  of the Dirac measure at the
origin. The distribution $ T$  satisfies the two conditions
$$ \langle T,f(v,v^*)\rangle =\langle T,f(tv,t^{-1}v^*)\rangle =|t|^n\langle T,f(tv,tv^*)\rangle .$$
therefore
$$ \langle T,f(v,t^2v^*)=|t|^{-n}\langle T,f(v,v^*)\rangle .$$

 Now $ T$   is recovered from $ R$  by the formula
$$ \langle T,f(v,v^*)=\int  _C\langle R,f(ce_n,^tc^{-1}v^*\rangle  dc=a\int  _Cf(ce_n,0)dc,\,\,\, f\in  {\mathcal  S}(\Omega  \times V^*).$$
Unless $ a=0$  this is not compatible with this last homogeneity condition.

Exactly in the same way one proves that $ T$  is 0 on  $ V \times
\Omega ^*$  where $ \Omega ^*$  is the open orbit $ x_1^*\ne 0$ of
$ C$  in $ V^*$. The same argument is valid for $ \widehat T$
(which is even equal to $ T$ \dots ).

 If $ n=1$  then $ T$  is obviously 0. If $ n\geq 2$  then there exists a distribution $ T'$  on
$$ \bigoplus_{1<j<n}{\mathbb   F}e_j\oplus{\mathbb   F}e_j^*$$
such that,
 $$ T=T'\otimes \delta _{x_n=0}\otimes dx_1\otimes \delta _{x_1^*=0}\otimes dx_n^*.$$
Let $ u$  be the isomorphism of $ V$  onto $ V^*$  given by $
u(e_j)=e_{n+1-j}^*$. Recall that it acts on $ {\mathfrak  g}
\times V \times V^*$  by $ (X,v,v^*)\mapsto
(^t(uXu^{-1}),^tu^{-1}v^*,uv)$. It belongs to $ \widetilde  C$ but
not to $ C$  so it must transform $ T$  into $ -T$.

The case $ n=1$  has just been settled. If $ n=2$   in the above
formula $ T'$  should be replaced by a constant. The constant must
be 0 if we want $ u(T)=-T$. If $ n>2$ let
$$V'=\Bigl (\oplus_1^{n-1}{\mathbb F}e_i\Bigr )/{\mathbb F}e_1$$
and let $X'$ be the nilpotent endomorphism of $V'$ defined by $X$.
We may consider $T'$ as a distribution on $V\oplus V^{'*}$ and one
easily checks that, with obvious notations, it transforms
according to the character $\chi$ of the the centralizer
$\widetilde C'$ of $X'$ in $\widetilde G'$. By induction $T'=0$,
hence $T=0$.

\hfill$\Box$\hfill

\section{\hglue -16pt. Reduction to the singular set: the orthogonal and unitary cases}

We now turn our attention to the unitary case. We keep the
notations of the introduction. In particular $W=V\oplus {\mathbb
D}e$ is a vector space over ${\mathbb  D}$ of dimension $n+1$ with
a non degenerate hermitian form $\langle.,.\rangle $ such that $e$
is orthogonal to $V$. The unitary group $G$ of $V$ is embedded
into the unitary group  $M$ of $W$.

Let $A$ be the set of all  bijective maps $u$ from $V$ to $V$ such that
$$ u(v_1+v_2)=u(v_1)+u(v_2),\,\,\, u(\lambda v)=\overline\lambda u(v),\,\,\, \langle u(v_1),u(v_2)\rangle =\overline{\langle v_1,v_2\rangle }.$$
An example of such a map is obtained by choosing a basis $e_1,\dots ,e_n$ of $V$ such that $\langle e_i,e_j\rangle \in {\mathbb  F}$ and defining
$$ u(\sum x_ie_i)=\sum \overline x_ie_i.$$
Any $u\in A$ is extended to $W$ by the rule $u(v+\lambda e)=u(v)+\overline\lambda e$ and we define an action on ${\mathrm GL}(W)$ by $m\mapsto um^{-1}u^{-1}$. The group $G$ acts on ${\mathrm GL}(W)$ by the adjoint action.

 Let $\widetilde G$ be the group of bijections of ${\mathrm GL}(W)$ onto itself generated by the actions of $G$  and $A$. It is a semi direct product of $G$ and $S_2$. We identify $G$ to a subgroup of $\widetilde G$ and $A$ to a subset. When a confusion is possible we denote the product in $\widetilde G$ with a $\times$.

  We define a character $\chi $ of $\widetilde G$ by $\chi (g)=1$ for $g\in G$ and $\chi (u)=-1$ for $u\in \widetilde G\setminus G$.
 Our overall goal is to prove that ${\mathcal  S}'(M)^{\widetilde G,\chi }=(0)$.

 Let $\widetilde G$ act on $G\times V$ as follows:
 $$ g(x,v)=(gxg^{-1},g(v)),\,\, u(x,v)=(ux^{-1}u^{-1},-u(v)),\quad g\in G,u\in A,x\in G,v\in V$$
Our first step is to replace $M$ by $G\times V$.
\begin{proposition}
 Suppose that for any $V$ and any hermitian form ${\mathcal  S}'(G\times V)^{\widetilde G,\chi }=(0)$, then ${\mathcal  S}'(M)^{\widetilde G,\chi }=(0)$.
\end{proposition}
{\it Proof.} We have in particular ${\mathcal  S}'(M\times
W)^{\widetilde M,\chi }=(0)$. Let $Y$ be the set of all $(m,w)$
such that $\langle w,w\rangle =\langle e,e\rangle $; it is a
closed subset, invariant under $\widetilde M$, hence ${\mathcal
S}'(Y)^{\widetilde M,\chi }=(0)$. By Witt's theorem $M$ acts
transitively on $\Gamma =\{w|\langle w,w\rangle =\langle
e,e\rangle \} $. We can apply Frobenius descent (Theorem 6-2) to the map
$(m,w)\mapsto w$ of $Y$ onto $\Gamma $. The centralizer of $e$ in
$\widetilde M$ is isomorphic to $\widetilde G$ acting as before on
the fiber  $M\times \{e\}$. We have a linear bijection between
${\mathcal  S}'(M)^{\widetilde G.\chi }$ and ${\mathcal
S}'(Y)^{\widetilde M,\chi }$; therefore ${\mathcal
S}'(M)^{\widetilde G.\chi }=(0)$.

\hfill$\Box$\hfill

The proof that ${\mathcal  S}'(G\times V)^{\widetilde G,\chi }=(0)$ is by induction on $n$. If ${\mathfrak  g}$ is the Lie algebra of $G$ we shall prove simultaneously that ${\mathcal  S}'({\mathfrak  g}\times V)^{\widetilde G,\chi }=(0)$. In this case $G$ acts on its Lie algebra by the adjoint action and for $u\in \widetilde G\setminus G$ one puts, for $X\in {\mathfrak  g},\,\,$ $u(X)=-uXu^{-1}$.

The case $n=0$ is trivial so we may assume that $n\geq 1$. If $T\in {\mathcal  S}'(G\times V)^{\widetilde G,\chi }$ in this section we will prove that the support of $T$ must be contained in the "singular set".

Let $Z$ (resp. {$\mathfrak  z$}) be the center of $G$ (resp. {$\mathfrak  g$}) and ${\mathcal  U}$ (resp. {$\mathcal  N$}) the (closed) set of all unipotent (resp. nilpotent) elements of $G$ (resp. {$\mathfrak  g$}).
\begin{lemma}
 If $T\in {\mathcal  S}'(G\times V)^{\widetilde G,\chi }$ (resp. $T\in {\mathcal  S}'({\mathfrak  g}\times V)^{\widetilde G,\chi }$) then the support of $T$ is contained in $Z{\mathcal  U}\times V$ (resp. ${\mathfrak  z}\times {\mathcal  N}\times V)$.
\end{lemma}
This is Harish-Chandra's descent. We first review some facts about
the centralizers of semi-simple elements, following \cite{SS}.

 Let $a\in G$, semi-simple; we want to describe its centralizer $M$ (resp. $\widetilde M$) in $G$ (resp. in $\widetilde G$) and to show that ${\mathcal  S}'(M\times
V)^{\widetilde M,\chi }=(0)$.

 View $a$ as a ${\mathbb  D}-$linear endomorphism of $V$ and call $P$ its minimal polynomial. Then, as $a$ is semi-simple, $P$ decomposes into irreducible factors $P=P_1\dots P_r$ two by two relatively prime. Let $V_i={\mathrm Ker } P_i(a)$ so that $V=\oplus V_i$. Any element $x$ which commutes with $a$ will satisfy $xV_i\subset V_i$ for each $i$.
For
$$ R(\xi )=d_0+\cdots +d_m\xi^m,\quad d_0d_m\ne  0$$
let
$$ R^*(\xi )=\overline{d_0}\xi ^m+\cdots +\overline{d_m}.$$
Then, from $aa^*=1$ we obtain, if $m$ is the degree of $P$
$$ \langle P(a)v,v'\rangle =\langle v,a^{-m}P^*(a)v'\rangle $$
(note that the constant term of $P$ can not be 0 because $a$ is
invertible). It follows that $P^*(a)=0$ so that $P^*$ is
proportional to $P$. Now $P^*=P_1^*\dots P_r^*$; hence there exists
a bijection $\tau $ from $\{1,2,\dots ,r\}$ onto itself such that
$P^*_i$ is proportional to $P_{\tau (i)}$. Let $m_i$ be the degree
of $P_i$. Then, for some non zero constant $c$
$$ 0=\langle P_i(a)v_i,v_j\rangle =\langle v_i,a^{-m_i}P_i^*(a)v_j\rangle =c\langle v_i,a^{-m_i}P_{\tau (i)}(a)v_j\rangle
,\quad v_i\in V_i,\,\, v_j\in V_j.$$
We have two possibilities.

{\bf Case 1:}$\,\, \tau (i)=i$. The space $V_i$ is orthogonal to
$V_j$ for $j\ne  i$; the restriction of the hermitian form to
$V_i$ is non degenerate. Let ${\mathbb  D}_i={\mathbb  D}[\xi
]/(P_i)$ and consider $V_i$ as a vector space over ${\mathbb D}_i$
through the action $(R(\xi ),v)\mapsto R(a)v$. As $a_{|V_i}$ is
invertible, $\xi $ is invertible modulo $(P_i)$; choose $\eta$
such that $\xi \eta=1$ modulo $(P_i)$. Let $\sigma _i$ be the
semi-linear involution of ${\mathbb  D}_i$, as an algebra over
${\mathbb  D}$:
$$ \sum d_j\xi^j\mapsto \sum \overline{d_j}\eta^j\quad{\mathrm modulo}\,\, (P_i)$$

Let ${\mathbb  F}_i$ be the subfield of fixed points for
$\sigma_i$. It is a finite extension of ${\mathbb  F}$, and
${\mathbb  D}_i$ is either a quadratic extension of ${\mathbb
F}_i$ or equal to ${\mathbb  F}_i$. There exists a ${\mathbb
D}-$linear form $\ell\ne  0$  on ${\mathbb  D}_i$ such that $\ell
(\sigma_i(d))=\sigma_i(\ell (d))$ for all $d\in {\mathbb  D }_i$.
Then any ${\mathbb  D}-$linear form $L$ on ${\mathbb  D}_i$ may be
written as $d\mapsto \ell (\lambda d)$ for some unique $\lambda
\in {\mathbb  D}_i$.

If $v,v'\in V_i$ then $d\mapsto \langle d(a)v,v'\rangle $ is ${\mathbb  D}-$linear map on ${\mathbb  D}_i$; hence there exists $S(v,v')\in {\mathbb  D}_i$ such that
$$ \langle d(a)v,v'\rangle =\ell (dS(v,v')).$$
One checks that $S$ is a non degenerate hermitian form on $V_i$ as
a vector space over ${\mathbb  D}_i$. Also a ${\mathbb  D}-$linear
map $x_i$ from $V_i$ into itself commutes with $a_i$ if and only
if it is ${\mathbb  D}_i$-linear and it is unitary with respect to
our original hermitian form if and only if it is unitary with
respect to $S$. So in this case we call $G_i$  the unitary group
of $S$. It does not depend upon the choice of $\ell$. As no
confusion may arise, for $\lambda \in {{\mathbb  D}_i}$ we define
$\overline\lambda =\sigma_i(\lambda )$.

We choose an ${\mathbb  F}_i-$linear map $u_i$ from $V_i$ onto
itself, such that $u_i(\lambda v)=\overline\lambda u(v)$ and
$S(u_i(v),u_i(v'))=\overline{S(v,v')}$. Then because of our
original choice of $\ell$ we also have $\langle
u_i(v),u_i(v')\rangle =\overline{\langle v,v'\rangle }$. Note that
$u(a_{V_i})^{-1}u^{-1}=a_{V_i}$.

{\bf Case 2.} Suppose now that $j=\tau (i)\ne  i$. Then $V_i\oplus V_j$ is orthogonal to $V_k$ for $k\ne  i,j$ and the restriction of the hermitian form to $V_i\oplus V_j$ is non degenerate, both $V_i$ and $V_j$ being totally isotropic subspaces. Choose an inverse $\eta$ of $\xi $ modulo $P_j$. Then for any $P\in {\mathbb  D}[\xi ]$
$$ \langle P(a)v_i,v_j\rangle =\langle v_i,\overline P(\eta (a))v_j\rangle ,\quad v_i\in V_i,\,\, v_j\in V_j$$
where $\overline P$ is the polynomial deduced from $P$ by changing
its coefficients into their conjugate. This defines a map, which
we call $\sigma_i$ from ${\mathbb  D}_i$ onto ${\mathbb  D}_j$. In
a similar way we have a map $\sigma _j$ which is the inverse of
$\sigma_i$. Then, for $\lambda \in {\mathbb  D}_i$ we have
$\langle \lambda v_i,v_j\rangle =\langle v_i,\sigma_i(\lambda
)v_j\rangle $.

 View $V_i$ as a vector space over ${\mathbb  D}_i$. The action
$$ (\lambda ,v_j)\mapsto
 \sigma_i(\lambda )v_j$$
defines a structure of ${\mathbb  D}_i$ vector space on $V_j$.
However note that for $\lambda \in {\mathbb  D}$ we have
$\sigma_i(\lambda )=\overline\lambda $ so that $\sigma_i(\lambda
)v_j$ may be different from $\lambda v_j$. To avoid confusion we
shall write, for $\lambda \in {\mathbb  D}_i$
$$\lambda v_i=\lambda *v_i\quad {\mathrm and}\quad \sigma_i(\lambda )v_j=\lambda *v_j.$$

As in the first case choose a non zero ${\mathbb  D}-$linear form $\ell$ on ${\mathbb  D}_i$. For $v_i\in V_i$ and $v_j\in V_j$ the map $\lambda \mapsto \langle \lambda *v_i,v_j\rangle $ is a ${\mathbb  D}-$linear form on ${\mathbb  D}_i$; hence there exists a unique element $S(v_i,v_j)\in {\mathbb  D}_i$ such that, for all $\lambda $
$$ \langle \lambda *v_i,v_j\rangle =\ell(\lambda S(v_i,v_j)).$$
The form $S$ is ${\mathbb  D}_i-$ bilinear and non degenerate so that we can view $V_j$ as the dual space over ${\mathbb  D}_i$ of the ${\mathbb  D}_i$ vector space $V_i$.

Let $(x_i,x_j)\in {\mathrm End}_{\mathbb  D}(V_i)\times {\mathrm End}_{\mathbb  D}(V_j)$. They commute with $(a_i,a_j)$ if and only if they are ${\mathbb  D}_i$-linear. The original hermitian form will be preserved, if and only if $S(x_iv_i,x_jv_j)=S(v_i,v_j)$ for all $v_i,v_j$. This means that $x_j$ is the inverse of the transpose of $x_i$. In this situation we define $G_i$ as the linear group of the ${\mathbb  D}_i-$vector space $V_i$.

Let $u_i$ be a ${\mathbb  D}_i-$linear bijection of $V_i$ onto $V_j$. Then $u_i(av_i)=a^{-1}u_i(v_i)$ and $u_i^{-1}(av_j)=a^{-1}u_i^{-1}(v_j)$.

Recall that $M$ is the centralizer of $a$ in $G$. Then $(M,V)$ decomposes as a "product", each "factor" being either of type $(G_i,V_i)$ with $G_i$ a unitary group
(case 1) or $(G_i,V_i\times V_j)$ with $G_i$ a general linear group (case 2). Gluing together the  $u_i$ (case 1) and the $(u_i,u_i^{-1})$ (case 2)
we get an element $u\in \widetilde G\setminus G$ such that $ua^{-1}u^{-1}=a$ which means that it belongs to the centralizer of $a$ in $\widetilde G$.
Finally if $\widetilde M$ is the centralizer of $a$ in $\widetilde G$ then $(\widetilde M,V)$ is  imbedded into  a product each "factor" being either
of type $(\widetilde G_i,V_i)$ with $G_i$ a unitary group (case 1) or $(\widetilde G_i,V_i\times V_j)$ with $G_i$ a general linear group (case 2).

If $a$ is not central then for each $i$ the dimension of $V_i$ is
strictly smaller than $n$ and from the result for the general
linear group and the inductive assumption in the  orthogonal or
unitary case we conclude that ${\mathcal  S}'(M\times
V)^{\widetilde M,\chi }=(0)$.

{\it Proof of Lemma 4.1.} in the group case. Consider the map
$g\mapsto P_g$ where $P_g$ is the characteristic polynomial of
$g$. It is a continuous map from $G$ into the set of polynomials
of degree at most $n$. Each non empty fiber ${\mathcal  F}$ is
stable under $ G$ but also under $\widetilde G\setminus G$.
Bernstein's localization principle tells us that it is enough to
prove that ${\mathcal  S}'({\mathcal  F}\times V)^{\widetilde
G,\chi }=(0)$.

Now it follows from \cite{SS} chapter IV that ${\mathcal  F}$
contains only a finite number of semi-simple orbits; in particular
the set of semi-simple elements ${\mathcal  F}_s$ in ${\mathcal
F}$ is closed. Let us use the multiplicative Jordan decomposition
into a product of a semi-simple and a unipotent element.  Consider
the map $\theta $ from ${\mathcal  F}\times V$ onto ${\mathcal
F}_s$ which associates to $(g,v)$ the semi-simple part $g_s$ of
$g$. This map is continuous (see the corresponding proof for $GL$)
and commutes with the action of $\widetilde G$. In ${\mathcal
F}_s$ each orbit $\gamma  $ is both open and closed therefore
$\theta ^{-1}(\gamma  )$ is open and closed and invariant under
$\widetilde G$. It is enough to prove that for each such orbit
${\mathcal  S}'(\theta ^{-1}(\gamma  ))^{\widetilde G,\chi }=(0)$.
By Frobenius descent (Theorem 6-2), if $a\in \gamma  $ and is not central, this
follows from the above considerations on the centralizer of such
an $a$ and the fact that $\theta ^{-1}(a)$ is a closed subset of
the centralizer of $a$ in $\widetilde G$, the product of the set
of unipotent element commuting with $a$ by $V$. Now $g_s$ is
central if and only if $g$ belongs to $Z{\mathcal  U}$, hence the
Lemma. For the Lie algebra the proof is similar, using the
additive Jordan decomposition.

\hfill$\Box$\hfill

Going back to the group if $a$ is central we see that it suffices to prove that ${\mathcal  S}'({\mathcal  U}\times V)^{\widetilde G,\chi }=(0)$ and similarly for the Lie algebra it is enough to prove that ${\mathcal  S}'({\mathcal  N}\times V)^{\widetilde G,\chi }=(0)$.

Now the exponential map (or the Cayley transform) is a
homeomorphism of ${\mathcal  N}$ onto ${\mathcal  U}$ commuting
with the action of $\widetilde G$. Therefore it is enough to
consider the Lie algebra case.

We now turn our attention to $V$. Let
$$ \Gamma =\{v\in V | \langle v,v\rangle =0\}.$$
\begin{proposition}
 If $T\in {\mathcal  S}'({\mathcal  N}\times V)^{\widetilde G,\chi }$ then the support of $T$ is contained in ${\mathcal  N}\times \Gamma $.
\end{proposition}
 {\it Proof.} Let
$$ \Gamma _t=\{v\in V\,|\, \langle v,v\rangle =0\}.$$
Each $\Gamma _t$ is stable by $\widetilde G$, hence, by
Bernstein's localization principle (Corollary 6-1), to prove that the support of
$T$ is contained in ${\mathcal  N}\times \Gamma _0$ it is enough
to prove that, for $t\ne  0$, $\,\, {\mathcal  S}'({\mathcal
N}\times \Gamma _t)^{\widetilde G,\chi }=(0)$.

By Witt's theorem the group $G$ acts transitively on $\Gamma _t$.
We can apply Frobenius descent to the projection from ${\mathcal
N}\times \Gamma _t$ onto $\Gamma _t$. Fix a point $v_0\in \Gamma
_t$. The fiber is ${\mathcal  N}\times \{v_0\}$. Let $\widetilde
G_1$ be the centralizer of $v_0$ in $\widetilde G$. We have to
show that ${\mathcal  S}'({\mathcal  N})^{\widetilde G_1,\chi
}=(0)$ and it is enough to prove that  ${\mathcal  S}'({\mathfrak
g})^{\widetilde G_1,\chi }=(0)$.

The vector $v_0$ is not isotropic so we have an orthogonal decomposition
$$ V={\mathbb  D}v_0\oplus V_1$$
with $V_1$ orthogonal to $v_0$. The restriction of the hermitian form to $V_1$ is non degenerate and $G_1$ is identified with the unitary group of this restriction, and   $\widetilde G_1$ is the expected semi-direct product with $S_2$. As a $\widetilde G_1-$module the Lie algebra ${\mathfrak  g}$ is isomorphic to a direct sum
$$ {\mathfrak  g\approx {\mathfrak  g}_1\oplus V_1\oplus W}$$
where ${\mathfrak  g}_1$ is the Lie algebra of $G_1$  and $W$ a vector space over ${\mathbb  F}$ of dimension 0 or 1 and on which the action of $\widetilde G_1$ is trivial. The action on ${\mathfrak  g}_1\oplus V_1$ is the usual one so that, by induction, we know that ${\mathcal  S}'({\mathfrak  g}_1\oplus V_1)^{\widetilde G_1,\chi }=(0)$. This readily implies that ${\mathcal  S}'({\mathfrak  g})^{\widetilde G_1,\chi }=(0)$.

\hfill$\Box$\hfill

Summarizing: we have to prove that ${\mathcal S}'({\mathcal N}\times\Gamma )^{\widetilde G,\chi}=(0)$.

\section{\hglue -16pt. End of the proof in the orthogonal and unitary cases}

We keep our general notations. We have to show that a distribution on ${\mathcal  N}\times \Gamma $ which is invariant under $G$ is invariant under $\widetilde G$. To some extent the proof will be similar to the one we gave for the general linear group.

 In particular we will use the fact that if $T$ is such a distribution then its partial Fourier transform on $V$ is also invariant under $G$. The Fourier transform on $V$ is defined using the bilinear form
 $$ (v_1,v_2)\mapsto \langle v_1,v_2\rangle +\langle v_2,v_1\rangle $$
 which is invariant under $\widetilde G$.

 For $v\in V$ put
 $$ \varphi   _v(x)=\langle x,v\rangle v,\quad x\in V.$$
 It is a rank one endomorphism of $V$ and $\langle \varphi   _v(x),y\rangle =\langle x,\varphi   _v(y)\rangle $.
 \begin{lemma}
 i)In the unitary case, for $\lambda \in {\mathbb  D}$ such that $\lambda =-\overline\lambda $ the map
 $$ \nu_\lambda : \quad (X,v)\mapsto (X+\lambda \varphi   _v,v)$$
 is a homeomorphism of $[{\mathfrak  g},{\mathfrak  g}]\times \Gamma $ onto itself which commutes with $\widetilde G$.\hfill\break
 ii) In the orthogonal case, for $\lambda \in {\mathbb  F}$ the map
 $$ \mu  _\lambda :\quad (X,v)\mapsto (X+\lambda X\varphi   _v+\lambda \varphi   _vX,v)$$
  is a homeomorphism of $[{\mathfrak  g},{\mathfrak  g}]\times \Gamma $ onto itself which commutes with $\widetilde G$.
 \end{lemma}
 The proof is a trivial verification.

 We now use the stratification of ${\mathcal  N}$. Let us first check that
 an adjoint orbit is stable not only by $G$ but by $\widetilde G$.

 Choose a basis $e_1,\dots ,e_n$ of $V$ such that $\langle e_i,e_j\rangle \in {\mathbb  F}$; this gives a conjugation $u: v=\sum x_ie_i\mapsto \overline v=\sum \overline{x_i}e_i$ on $V$. If $A$ is any endomorphism of $V$ then $\overline A$ is the endomorphism $v\mapsto \overline{A(\overline v)}$. The conjugation $u$ is an element of $\widetilde G\setminus G$ and, as such, it acts on ${\mathfrak  g}\times V$ by $(X,v)\mapsto (-uXu^{-1},-u(v))=(-\overline X,-\overline v)$. In \cite{MVW} Chapter 4 Proposition 1-2
it is shown that for $X\in {\mathfrak  g}$ there exists an
${\mathbb  F}-$linear automorphism $a$ of $V$ such that $\langle
a(x),a(y)\rangle =\overline{\langle x,y\rangle }$ (this implies
that $a(\lambda x)=\overline{\lambda }x$) and such that
$aXa^{-1}=-X$. Then $g=ua\in G$ and $gXg^{-1}=-\overline X$ so
that $-\overline X$ belongs to the adjoint orbit of $X$. Note that
$a\in \widetilde G\setminus G$ and as such acts as
$a(X,v)=(X,-a(v))$; it is an element of the centralizer of $X$ in
$\widetilde G\setminus G$.

{\textbf {Remark.}} We need to check this only for nilpotent
orbits and this will be done later in an   explicit way, using the
canonical form of nilpotent matrices.

  Let ${\mathcal  N}_i$ be the union of all nilpotent orbits of dimension at most $i$. We shall prove, by descending induction on $i$, that the support of a distribution $T\in {\mathcal  S}'({\mathcal  N}\times \Gamma )^{\widetilde G,\chi }$ must be contained in ${\mathcal  N}_i\times \Gamma $.

  So now assume that $i\geq 0$ and that we already know that the support of any $T\in {\mathcal  S}'({\mathcal  N}\times \Gamma )^{\widetilde G,\chi }$ must be contained in ${\mathcal  N}_i\times \Gamma $. Let ${\mathcal  O}$ be a nilpotent orbit of dimension $i$; we have to show that the restriction of $T$ to ${\mathcal  O}$ is 0.

In the unitary case fix $\lambda \in {\mathbb  D}$ such that
$\lambda =-\overline\lambda $ and consider,  for every $t\in
{\mathbb F}$ the homeomorphism $\nu_{t\lambda }$; the image of $T$
belongs to ${\mathcal  S}'({\mathcal  N}\times \Gamma
)^{\widetilde G,\chi }$ so that the image of the support of $T$
must  be contained in ${\mathcal  N}_i\times \Gamma $. If $(X,v)$
belongs to this support this means that $X+t\lambda \varphi _v\in
{\mathcal N}_i$.

If $i=0$ so that $\mathcal {N}_i=\{0\}$ this implies that  $v=0$ so that $T$ must be a multiple of the Dirac measure at the point $(0,0)$ and hence is invariant under $\widetilde G$ so must be 0.

If $i>0$ and $X\in {\mathcal  O}$ then as ${\mathcal  O}$ is open in ${\mathcal  N}_i$, we get that, at least for $|t|$ small enough, $X+t\lambda \varphi   _v\in {\mathcal  O}$ and therefore $\lambda \varphi   _v$ belongs to the tangent space $\mathrm{Im}\,\mathrm{ad}   (X)$ of ${\mathcal  O}$ at the point $X$. Define
$$ Q(X)=\{v\in V|\varphi   _v\in \mathrm {Im}\,\mathrm {ad }  (X)\},\quad X\in {\mathcal  N},\,\,\, {\mathrm {(unitary \,\, case)}}.$$
Then we know that the support of the restriction of $T$ to ${\mathcal  O}$ is contained in
$$ \{(X,v)|X\in {\mathcal  O,\, v\in Q(X)}\}$$
and the same is true for the partial Fourier transform of $T$ on $V$.

In the orthogonal case for $i=0$, the distribution $T$ is the product of the Dirac measure at the origin of $\mathfrak{g}$ by  a distribution $T'$ on $V$. The distribution $T'$ is invariant under $G$ but the image of $\widetilde G$ in $\mathrm{End} (V)$ is the same as the image of $G$ so that $T'$ is invariant under $\widetilde G$ hence must be 0.

If $i>0$ we proceed  as in the unitary case, using $\mu  _\lambda $. We define
$$ Q(X)=\{v\in V|X\varphi   _v+\varphi   _vX\in {\mathrm Im}\,{\mathrm ad }  (X)\},\quad X\in {\mathcal  N},\,\,\, {\mathrm {(orthogonal \,\, case)}}$$
and we have the same conclusion.

In both cases, for $i>0$, fix $X\in {\mathcal  O}$. We use
Frobenius descent for the projection  map  $(Y,v)\mapsto Y$ of
${\mathcal  O}\times V$ onto ${\mathcal  O}$. Let $C$ (resp.
$\widetilde C$) be the stabilizer of $X$ in $G$ (resp. $\widetilde
G$). We have a linear bijection of ${\mathcal  S}'({\mathcal
O}\times \Gamma )^{\widetilde G,\chi }$ onto ${\mathcal
S}'(V)^{\widetilde C,\chi }$.
\begin{lemma}
  Let $T\in {\mathcal  S}'(V)^{\widetilde C,\chi }$. If $T$ and its Fourier transform are supported in $Q(X)$ then $T=0$.
\end{lemma}
Let us say that a nilpotent element $X$ is nice if the above Lemma is true.

Suppose that we have a direct sum decomposition $V=V_1\oplus V_2$
such that $V_1$ and $V_2$ are orthogonal. By restriction we get
non degenerate hermitian forms $\langle .,.\rangle _i$ on $V_i$.
We call $G_i$ the unitary group of $\langle .,.\rangle _i$,
${\mathfrak  g}_i$ its Lie algebra and so on. Suppose that
$X(V_i)\subset V_i$ so that $X_i=X_{|V_i}$ is a nilpotent element
of ${\mathfrak  g}_i$.
\begin{lemma}
 If $X_1$ and $X_2$ are nice so is $X$.
\end{lemma}
{\it Proof of Lemma 5.3. }We claim that $Q(X)\subset Q(X_1)\times Q(X_2)$. Indeed if
$$ A=\pmatrix{A_{1,1}&A_{1,2}\cr A_{2,1}&A_{2,2}\cr}\in {\mathfrak  g}$$
then from
$$ \langle A\pmatrix{x_1\cr x_2\cr},\pmatrix{y_1\cr y_2\cr}\rangle +\langle \pmatrix{x_1\cr x_2\cr},A\pmatrix{y_1\cr y_2\cr}\rangle =0$$
we get in particular
$$ \langle A_{i,i}x_i,y_i\rangle +\langle x_i,A_{i,i}y_i\rangle =0$$
so that $A_{i,i}\in {\mathfrak  g}_i$. Note that
$$ [X,A]=\pmatrix{[X_1,A_{1,1}]&*\cr *&[X_2,A_{2,2}]\cr}.$$

If $v_i\in V_i$ and $v_j\in V_j$ we define $\varphi   _{v_i,v_j}: V_i\mapsto V_j$ by $\varphi   _{v_i,v_j}(x_i)=\langle x_i,v_i\rangle v_j$. Then, for $v=v_1+v_2$
$$ \varphi   _{v}=\pmatrix{\varphi   _{v_1,v_1}&\varphi   _{v_2,v_1}\cr \varphi   _{v_1,v_2}&\varphi   _{v_1,v_2}\cr}.$$
Therefore if, for $A\in {\mathfrak  g}$ we have $\varphi   _v=[X,A]$ then $\varphi   _{v_i,v_i}=[X_i,A_{i,i}]$. This proves the assertion for the unitary case. The orthogonal case is similar.

The  end of the proof  is the same as the end of the proof of Lemma 3-2.

\hfill$\Box$\hfill

Now in both orthogonal and unitary cases nilpotent elements have normal forms which are orthogonal direct sums of "simple" nilpotent matrices. This is precisely described in \cite{SS} IV 2-19 page 259. By the above Lemma it is enough to prove that each "simple" matrix is nice.

{\bf Unitary case.} There is only one type to consider.  There
exists a basis $e_1,\dots ,e_n$ of $V$ such that $Xe_1=0$ and
$Xe_i=e_{i-1}, \, i\geq 2$. The hermitian form is given by
$$ \langle e_i,e_j\rangle =0\, \, {\mathrm if}\, i+j\ne  n+1,\quad \langle e_i,e_{n+1-i}\rangle =(-1)^{n-i}\alpha$$
with $\alpha\ne  0$. Note that $\overline\alpha =(-1)^{n-1}\alpha
$. Suppose that $v\in Q(X)$; for some $A\in {\mathfrak  g}$ we
have $\lambda \varphi   _v=XA-AX$. For any integer $p\geq 0$
$${\mathrm Tr  }  (\lambda \varphi   _vX^p)={\mathrm Tr} (XAX^{p}-AX^{p+1})=0.$$
Now
${\mathrm Tr} (\varphi   _vX^p)=\langle X^pv,v\rangle $
Let $v=\sum x_ie_i$. Hence
$$ \langle X^pv,v\rangle =\sum _1^{n-p}x_{i+p}\langle e_i,v\rangle =\sum _1^{n-p}(-1)^{n-i}\alpha x_{i+p}\overline x_{n+1-i}=0.$$
For $p=n-1$ this gives $x_n\overline x_n=0$. For $p=n-2$ we get
nothing new but for $p=n-3$ we obtain $x_{n-1}=0$. Going on, by an
easy induction, we conclude that $x_i=0$ if $i\geq (n+1)/2$.

If $n=2p+1$ is odd put $V_1=\oplus _1^p{\mathbb  D}e_i$,$\,\, V_0={\mathbb  D}e_{p+1}$ and $V_2=\oplus _{p+2}^{2p+1}{\mathbb  D}e_i$. If $n=2p$ is even put $V_1=\oplus _1^p{\mathbb  D}e_i$, $V_0=(0)$ and $V_2=\oplus _{p+1}^{2p}{\mathbb  D}e_i$. In both cases we have $V=V_1\oplus V_0\oplus V_2$. We use the notation $v=v_2+v_0+v_1$

The distribution $T$ is supported by $V_1$. Call $\delta  _i$ the
Dirac measure at $0$ on $V_i$. Then we may write $T=U\otimes
\delta _0\otimes \delta  _2$ with $U\in {\mathcal  S'(V_1)}$. The
same thing must be true of the Fourier transform of $T$. Note that
$\widehat U$ is a distribution on $V_2$, that $\widehat\delta _2$
is a Haar measure $dv_1$ on $V_1$ and that, for $n$ odd
$\widehat\delta _0$ is a Haar measure $dv_0$ on $V_0$. So we have
$\widehat T= dv_1\otimes \widehat U$ if $n$ is even and $\widehat
T=dv_1\otimes dv_0\otimes \widehat U$ if $n$ is odd. In the odd
case this forces $T=0$. In the even case, up to a scalar multiple
the only possiblity is $T=dv_1\otimes \delta  _2$.

Let $$ a: \sum x_ie_i\mapsto \sum (-1)^i\overline x_ie_i.$$ Then
$a\in \widetilde G\setminus G$. It acts on ${\mathfrak  g}$ by
$Y\mapsto -aYa^{-1}$ and in particular  $-aXa^{-1}=X$ so that
$a\in \widetilde C\setminus C$. The action on $V$ is given by
$v\mapsto -a(v)$. It is an involution. The subspace $V_1$ is
invariant and so $dv_1$ is invariant. This implies that $T$ is
invariant under $\widetilde C$ so it must be 0.

{\bf Orthogonal case.} There are two different types of "simple"
nilpotent matrices.

{\bf The first type} is the same as the unitary case, with
$\alpha=1$ and thus $n$ odd but now our condition is that
$X\varphi _v+\varphi   _vX=[X,A]$
 for some $A\in {\mathfrak  g}$. As before this implies that ${\mathrm Tr} (\varphi   _vX^q)=0$ but only for $q\geq 1$. Put $n=2p+1$; we get $x_j=0$ for $j>p+1$.
Decompose $V$ as before: $V=V_1\oplus V_0\oplus V_2$. Our
distribution $T$ is supported by the subspace $v_2=0$ so we write
it $T=U\otimes \delta  _2$ with $U\in {\mathcal  S}'(V_1\oplus
V_0)$. This is also true  for the distribution $\widehat T$ so we
must have $U=dv_1\otimes R$ with $R$ a distribution on $V_0$.
Finally $T=dv_1\otimes R\otimes \delta  _2$. Now $-{\mathrm Id}\in
C$ and $T$ is invariant under $C$ so that $R$ must be an even
distribution. On the other end the endomorphism $a$ of $V$ defined
by $a(e_i)=(-1)^{i-p-1}e_i$ belongs to $C$and $aXa^{-1}=-X$ and
$u: (X,v)\mapsto (-X,-v)$ belongs to $\widetilde G\setminus G$.
The product $a*u$ of $a$ and $u$ in $\widetilde G$ belongs to
$\widetilde C \setminus C$. Clearly $T$ is invariant under $a*u$
so that $T$ is invariant under $\widetilde C$ so it must be 0.

{\bf The second type} is as follows. We have $n=2m$, an even integer and a decomposition $V=E\oplus F$ with both $E$ and $F$ of dimension $m$. We have a basis $e_1,\dots ,e_m$ of $E$ and a basis $f_1,\dots ,f_m$ of $F$ such that
$$ \langle e_i,e_j\rangle =\langle f_i,f_j\rangle =0$$
and
$$ \langle e_i,f_j\rangle =0 \,{\mathrm if}\, i+j\ne  m+1\quad {\mathrm and}\quad \langle e_i,f_{m+1-i}\rangle =(-1)^{m-i}.$$
Finally $X$ is  such that $Xe_i=e_{i-1},\,\, Xf_i=f_{i-1}$.

Let $\xi $ be the matrix of the restriction of $X$ to $E$ or to $F$. Write an element $A\in {\mathfrak  g}$ as 2 by 2 matrix $A=(a_{i,j})$. Then
$$[X,A]=\pmatrix{[\xi ,a_{1,1}]&[\xi ,a_{1,2}]\cr [\xi ,a_{2,1}]& [\xi ,a_{2,2}]\cr}.$$
 Suppose that $v\in Q(X)$ and let
$$ v=e+f\,\,\, {\mathrm with}\,\,\,e=\sum x_ie_i,\,\,\,f=\sum y_if_i.$$
We get
$$ X\varphi   _v+\varphi   _vX=\pmatrix{\xi \varphi   _{f,e}+\varphi   _{f,e}\xi &\xi \varphi   _{f,f}+\varphi   _{f,f}\xi \cr\xi \varphi   _{e,e}+\varphi   _{e,e}\xi &\xi \varphi   _{e,f}+\varphi   _{e,f}\xi \cr}$$
where, for example $\varphi   _{e,e}$ is the map $f'\mapsto \langle f',e\rangle e$ from $F$ into $E$. Thus, for some $A$,
$$\xi \varphi   _{e,e}+ \varphi   _{e,e}\xi =\xi a_{2,1}-a_{2,1}\xi$$
In this formula, using the basis $(e_i),\, (f_i)$ replace all the
maps by their matrices.

 Then, as before, we have ${\mathrm Tr} (\varphi   _{e,e}\xi ^q)=0$ for $1\leq q\leq   m-1$. If $e'=\sum x_if_i$ (the $x_i$ are the coordinates of $e$),
 then ${\mathrm Tr} (\xi^q\varphi   _{e,e})$ is $\langle\xi ^qe,e'\rangle $. Thus, as in the other cases, we have $x_j=0$ for $j>n/2$ if $n$ is even and $j>(n+1)/2$ if $n$ is odd. The same thing is true for the $y_i$.

 If $n=2p$ is even, let $V_1=\oplus _{i\leq   p}({\mathbb  F}e_i\oplus {\mathbb
F}f_i)$ and $V_2=\oplus _{i>p}({\mathbb  F}e_i\oplus {\mathbb
F}f_i)$; write $v=v_1+v_2$ the corresponding decomposition of an
arbitrary element of $V$. Let $\delta _2$ be the Dirac measure at
the origin in $V_2$ and $dv_1$ a Haar measure on $V_1$. Then, as
in the unitary case, using the Fourier transform, we see that the
distribution $T$ must be a multiple of $dv_1\otimes \delta_2$.

 The endomorphism $a$ of $V$ defined by $a(e_i)=(-1)^ie_i$ and $a(f_i)=(-1)^{i+1}f_i$ belongs to $G$ and $aXa^{-1}=-X$. The map $u: (Y,v)\mapsto (-Y,-v)$ belongs to $\widetilde G\setminus G$ so that the product $a\times u$ in $\widetilde G$ belongs to $\widetilde C\setminus C$. It clearly leaves $T$ invariant so that $T=0$.

 Finally if $n=2p+1$ is odd we put $V_1=\oplus _{i\leq   p}({\mathbb  F}e_i\oplus {\mathbb  F}f_i)$, $V_0={\mathbb  F}e_{p+1}\oplus {\mathbb  F}f_{p+1}$,
  $V_2=\oplus _{i\geq p+2}({\mathbb  F}e_i\oplus {\mathbb  F}f_i)$.
 As in the unitary case we find that $T=dv_1\otimes R\otimes \delta  _2$ with $R$ a distribution on $V_0$. As $-{\mathrm Id}\in C$ we see that $R$ must be even.
 Then again, define $a\in G$ by $a(e_i)=(-1)^ie_i$ and $a(f_i)=(-1)^if_i$ and consider $a*u$ with $u(Y,v)=(-Y,-v)$.
 As before $a*u\in \widetilde C\setminus C$ and leaves $T$ invariant so we have to take $T=0$.

 \hfill$\Box$\hfill

\section{Appendix}
We shall state two theorems which are systematically used in our
proof.

If $X$ is a Hausdorff totally disconnected  locally compact
topological space (lctd space in short) we denote by ${\mathcal
S}(X)$ the vector space of locally constant applications with
compact support of $X$ into the field of complex numbers ${\mathbb
C}$. The dual space ${\mathcal S}'(X)$ of ${\mathcal S}(X)$ is the
space of distributions on $X$. All the lctd spaces we introduce
are countable at infinity.

If an lctd topological group  $G$ acts continuously on a lctd
space $X$ then it acts on ${{\mathcal S}}(X)$ by \[
(gf)(x)=f(g^{-1}x)\] and on distributions by
\[ (gT)(f)=T(g^{-1}f)\]
The space of invariant distributions is denoted by ${\mathcal
S}'(X)^G$. More generally, if $\chi$ is a character of $G$ we
denote by ${\mathcal S}'(X)^{G,\chi}$ the space of distributions
$T$ which transform according to $\chi$ that is to say
$T(f(g^{-1}x))=\chi (g)T(f)$.

The following result is due to Bernstein \cite{Ber}, section 1.4.
\begin{theorem}[Localization principle] Let  $q:Z \to T$ be a continuous map between two topological spaces of type lctd. Denote $Z_t:=
q^{-1}(t)$. Consider $\mathcal{S}'(Z)$ as $\mathcal{S}(T)$-module.
Let $M$ be a closed subspace of $\mathcal{S}'(Z)$ which is an
$\mathcal{S}(T)$-submodule. Then $M=\overline{\bigoplus_{t \in T}
(M \cap \mathcal{S}'(Z_t))}$.
\end{theorem}

\begin{corollary} Let $q:Z \to T$ be a continuous map between topological spaces of type lctd. Let an lctd group
$H$ act on  $Z$ preserving the fibers of $q$. Let $\mu$ be a
character of $H$. Suppose that for any $t\in T$,
$\mathcal{S}'(q^{-1}(t))^{H,\mu}=0$. Then
$\mathcal{S}'(Z)^{H,\mu}=0$.
\end{corollary}

The second theorem is a variant of Frobenius reciprocity.

\begin{theorem}[Frobenius descent]
Let a unimodular lctd  topological group $H$ act transitively on
an  lctd topological space $Z$. Let $\varphi:E \to Z$ be an
$H$-equivariant map of lctd topological spaces. Let $x\in Z$.
Suppose that its stabilizer $Stab_H(x)$ is unimodular. Let $W$ be
the fiber of $x$.
Let $\chi$ be a character of $H$. Then\\
(i) There exists a canonical isomorphism $Fr:
\mathcal{S}'(E)^{H,\chi}
\to \mathcal{S}'(W)^{Stab_H(x),\chi}$.\\
(ii) For any distribution $\xi \in \mathcal{S}'(E)^{H,\chi}$,
$\mathrm{Supp}(Fr(\xi))=\mathrm{Supp}(\xi)\cap W$.\\
(iii) Frobenius descent commutes with Fourier transform.

Namely, let $W$ be a finite dimensional linear space over $\mathbb
F$ with a nondegenerate bilinear form $B$. Let $H$ act on $W$
linearly preserving $B$.

 Then
for any $\xi \in \mathcal{S}'(Z\times W)^{H,\chi}$, we have
$\mathcal{F}_{B}(\mathrm{Fr}(\xi))=\mathrm{Fr}(\mathcal{F}_{B}(\xi))$
where $\mathrm{Fr}$ is taken with respect to the projection $Z
\times W \to Z$.
\end{theorem}

\vskip 1cm

\noindent A. Aizenbud \\
Faculty of Mathematics and Computer Science, Weizmann Institute of
Science POB 26, Rehovot 76100, ISRAEL. \\
E-mail: aizenr@yahoo.com

$ $\\
D. Gourevitch\\
Faculty of Mathematics and Computer Science, Weizmann Institute of
Science POB 26, Rehovot 76100, ISRAEL.\\
E-mail: guredim@yahoo.com.

$ $\\
S. Rallis\hfill\break Department of Mathematics Ohio State
University,\hfill\break COLUMBUS, OH 43210\hfill\break
E-mail:haar@math.ohio-state.edu

$ $\\
G. Schiffmann\hfill\break Institut de Recherche Math\'ematique
Avanc\'ee,\hfill\break Universit\'e Louis-Pasteur et CNRS, 7, rue
Ren\'e-Descartes,\hfill\break 67084 Strasbourg Cedex,
France.\hfill\break e-mail: schiffma@math.u-strasbg.fr

 \end{document}